\documentclass[11pt]{article}
 \usepackage[a4paper, total={6.7in, 8.7in}]{geometry}
\usepackage{setspace}
\usepackage[utf8]{inputenc}
\usepackage[utf8]{inputenc}
\usepackage{amssymb}
\usepackage{amsmath}
\usepackage{amsthm}
\usepackage{hyperref}
\usepackage{breqn}
\usepackage[table,xcdraw]{xcolor}
\usepackage{xcolor}
\usepackage{mathtools}
\usepackage{bm}
\usepackage{enumitem}
\usepackage{mathrsfs}
\usepackage{graphicx}
\usepackage{caption}
\usepackage{url}
\usepackage{euscript}
\usepackage{placeins}
\usepackage{mathdots}
\usepackage{multirow}
\usepackage{arydshln} 
\usepackage[english]{babel}
\usepackage{array}
\usepackage{subcaption}
\usepackage{booktabs}
\usepackage[numbers]{natbib}
\usepackage{cases}

\usepackage{graphics}
\usepackage{graphicx}
\usepackage{algorithm}
\usepackage{algpseudocode}
\numberwithin{equation}{section}
\newcommand{\dm}[1]{{\displaystyle{#1}}}

\newcommand{\iu}{{i\mkern1mu}}
\newtheorem{theorem}{\bf Theorem}[section]

\newtheorem{corollary}{Corollary}[section]
\newtheorem{proposition}{Proposition}[section]
\newtheorem{lemma}{Lemma}[section]
\newtheorem{remark}{Remark}[section]
\newtheorem{method}{Method}[section]
\theoremstyle{remark}
\newtheorem{exam}{\bf Example}

\def \R{{\mathbb R}}
\def\bmatrix#1{\left[\begin{matrix}
		#1
	\end{matrix}\right]}
\def \diag{\mathrm{diag}}
\def \B{\mathcal B}
\def \P{\mathcal P}

\def \om{\bm \omega}
\def \z{\bm z}
\def \w{\bm w}

\def \y{\bm y}
\def \0{\bm 0}
\def \GSS{\mathrm{GSS}}

\def\bmatrix#1{\left[ \begin{matrix} #1 \end{matrix} \right]}

\def \R{{\mathbb R}}

\title{A Class of Generalized Shift-Splitting Preconditioners for Double Saddle Point Problems}
\author{Sk. Safique Ahmad\footnotemark[2] \footnotemark[1] \and Pinki Khatun \footnotemark[2] }
\date{}
\begin{document}
\onehalfspacing
	\maketitle	
	\begin{abstract}
	In this paper, we propose a generalized shift-splitting (GSS) preconditioner, along with its two relaxed variants to solve the double saddle point problem (DSPP). The convergence of the associated GSS iterative method is analyzed, and sufficient conditions for its convergence are established. Spectral analyses are performed to derive sharp bounds for the eigenvalues of the preconditioned matrices. Numerical experiments based on examples arising from the PDE-constrained optimization problem and the leaky lid driven cavity problem demonstrate the effectiveness and robustness of the proposed preconditioners compared with existing state-of-the-art preconditioners.
	\end{abstract}
	
	

	\footnotetext[2]{
		Department of Mathematics, Indian Institute of Technology Indore, Khandwa Road, Indore 453552, Madhya Pradesh, India.  \texttt{Email:\,safique@iiti.ac.in}, \texttt{pinki996.pk@gmail.com}}
	\footnotetext[1]{Corresponding author.
		}

\section{Introduction}
 Suppose $n,m$, and $l$ are given positive integers with $n\geq m\geq l.$ Then, we consider the double saddle point problem (DSPP) in the following form \cite{Susanne2023}:
\begin{equation}\label{SPP1}
    \mathcal{B}{\bm u}:=\bmatrix{A &\0&B^T\\ \0&D&C\\-B&-C^T&\0}\bmatrix{x\\y\\z}=\bmatrix{p\\q\\r}=:{\bm b},
\end{equation}
where $A\in \R^{n\times n}$, $B\in \R^{m\times n}$ has full row rank, $C\in \R^{l\times m}$ and $D\in \R^{l\times l}$ . Further, $p\in \R^n,$ $q\in \R^l$ and $r\in \R^{m}$ are known vectors and   $x\in \R^n,$ $y\in \R^l$ and $z\in \R^{m}$ are unknown vectors to be determined. Here, $(\cdot)^T$ and $\0$ represent the transpose of a matrix and zero matrix of appropriate size, respectively.  {In this paper, we consider that the matrices \( A \) and \( D \) can be both symmetric or nonsymmetric.}

The DSPP \eqref{SPP1}  is frequently encountered in a wide range of scientific and computational disciplines. Notable areas of application include quadratic programming problems \cite{OPTMIZATION1}, computational fluid dynamics \cite{CFD2005},  equality constrained indefinite least squares problems \cite{ILSE2003}, PDE-constrained optimization problem \cite{PDE-constrained2010}, and so on.

Owing to the broad applicability of the DSPP \eqref{SPP1}, this paper primarily concentrates on its numerical solution. Nevertheless, the large and sparse nature of the saddle point matrix $\B$ makes the iterative methods more favorable than the direct approaches \cite{YSAAD}. In particular, Krylov subspace methods, such as the generalized minimal residual (GMRES) method \cite{gmres}, have emerged prominently as one of the most effective and widely used iterative solvers. However, the ill-conditioned nature and poor spectral properties of the coefficient matrix $\B$ lead to a significant slowdown of the Krylov subspace method. This challenge underscores the need for the development of novel, effective, and robust preconditioners to enhance the performance of Krylov subspace methods.

Assuming $\mathbf{A}=\bmatrix{A& \0 \\ \0 & D},$ $\mathbf{B}=[B ~~ C^T],$ ${\bf f}=[p^T, q^T]^T,$ ${\bf g}=r,$ ${\bf x}=[x^T, y^T]^T$ and ${\bf y}=z,$ the DSPP can be transformed into the two-by-two block saddle point problems \cite{Benzi2005}.
A surge amount of work has been done in recent years to solve two-by-two block saddle point problems, for instance,  Uzawa iterative methods \cite{ZZBaiUzawa, MUzawa2014, Uzawa2014}, shift-splitting (SS) methods \cite{CaoSS2014, CaoSS2017}, successive overrelaxation methods \cite{GSOR2005, SOR2001}, Hermitian and skew-Hermitian methods \cite{ HSS2019, HSS}, and so on. However, these methods do not exploit the three-by-three block structure of $\B.$

To leverage the full block structure of $\B,$ various preconditioners have been studied in the literature to solve the DSPP \eqref{SPP1}. When $D=\0$ in \eqref{SPP1}, block diagonal (BD) and tridiagonal type preconditioners \cite{BDMaryam, HuangNA}, SS type preconditioners \cite{Pinki_PESS, CAOSS19}, Uzawa methods \cite{DUM2019, NaHuangVUZ}, etc. have been explored. When $D\neq \0,$ BD preconditioners for the DSPP \eqref{SPP1} have been investigated in \cite{Susanne2023}.   
By splitting the coefficient matrix $\B$ as $\B=\B_1+\B_2,$  where  \begin{eqnarray}
   \B_1=\bmatrix{A & \0 & B^T\\ \0 & \0 & \0\\ -B & \0 & \0} ~ \text{and} ~~   \B_2=\bmatrix{\0 & \0 & \0 \\ \0 & D& C\\ \0 & -C^T & \0},
\end{eqnarray}
\citet{BenziDS2011} proposed the dimensional spitting (DS) preconditioner $\mathcal{P}_{\mathrm{DS}}$ and a relaxed dimensional factorization (RDF) preconditioner $\mathcal{P}_{\mathrm{RDF}}$ \cite{BenziRDF2011}. These are given as follows: 
\begin{align}
&\mathcal{P}_{\mathrm{DS}}=\frac{1}{\alpha} \bmatrix{\alpha I+A & \0 & B^T \\ \0 & \alpha I & \0 \\ -B & \0 & \alpha I} \bmatrix{ \alpha I & \0 & \0 \\ \0 & \alpha I+ D & C\\ \0 & -C^T & \alpha I},\\
&\mathcal{P}_{\mathrm{RDF}}=\frac{1}{\alpha} \bmatrix{A & \0 & B^T \\ \0 & \alpha I & \0 \\ -B & \0 & \alpha I} \bmatrix{ \alpha I & \0 & \0 \\ \0 &  D & C\\ \0 & -C^T & \alpha I},
\end{align}
where $I$ stands for the identity matrix of appropriate size. For more on DS-based preconditioners, refer to \cite{SDS2019, MDS2024}. 

For the DSPP \eqref{SPP1} arising from PDE-constrained optimization problem, \citet{PDE-constrained2010} introduced BD preconditioner $\mathcal{P}_{\text{D}}$ and constrained preconditioner $\mathcal{P}_{\text{C}}.$ Later, a block triangular (BT) preconditioner is developed in \cite{ReesBT2010}. In \cite{BCD2014}, the authors introduced two types of preconditioners for solving the DSPP: a block-counter-diagonal preconditioner, denoted as $\mathcal{P}_{\mathrm{BCD}}$, and a block-counter-tridiagonal preconditioner, denoted as $\mathcal{P}_{\mathrm{BCT}}$. By using different approximations of Schur complement, preconditioners in BD and BT formats are constructed in \cite{Pearson2012, Pearson2, Pearson3}. For more research on solving DSPP in the context of PDE-constrained optimization problems, refer to \cite{SSPDE2024, ChangFeng2018, Salkuyeh2020}.

For non-Hermitian linear systems, \citet{BaiSS} introduced the SS preconditioner by applying the matrix SS technique to the coefficient matrix. Building on this concept, several SS-type preconditioners have been studied for two-by-two and three-by-three block saddle point problems  {with nonsymmetric coefficient matrix}; see, for example, \cite{CAOSS19, CaoSS2014, CaoSS2017, GSS2015, Pinki_PESS}. However, SS-type preconditioners remain unexplored for DSPP, despite their potential for high efficiency. Motivated by this, in this paper, we introduce a generalized SS (GSS) iterative method and the associated GSS preconditioner for solving the DSPP by implementing the concept of SS technique to $\B$. Moreover, we derive sufficient criteria for the convergence of the GSS iterative method. The main contributions of the paper are summarized as follows:
  \begin{itemize}
     \item A novel GSS iterative method and corresponding GSS preconditioner are introduced by implementing the SS approach for the coefficient matrix $\B$ to solve DSPP \eqref{SPP1}.
     \item Convergence analysis of the proposed GSS iterative method are carried out, yielding sufficient conditions for its convergence.
      \item To further enhance the effectiveness of the GSS preconditioner, two relaxed variants, termed RGSS-I and RGSS-II are proposed, and the spectral bounds of the RGSS-I and RGSS-II preconditioned matrices are thoroughly investigated.
      \item Finally, numerical experiments are conducted for the DSPP arising from the PDE-constrained optimization problem and  {the leaky lid driven cavity problem} to demonstrate the effectiveness of the proposed preconditioners.
  \end{itemize}
 

The outline of the rest of the paper is as follows.  Section \ref{sec2} investigates the solvability conditions of the DSPP \eqref{SPP1} and properties of the coefficient matrix $\B.$ The GSS iterative method and associated preconditioner are proposed in Section \ref{sec:GSS}. Section \ref{SEC:convergence} investigates the convergence criteria for the proposed GSS iterative method. Two relaxed variants of the proposed GSS preconditioner are presented in Section \ref{sec:RGSS}, and the spectral analyses of the corresponding preconditioned matrices are performed. Section \ref{sec:parameter} deals with the parameter selection strategies of the proposed preconditioners. Experimental results, analyses and discussions of the proposed and existing preconditioners are discussed in Section \ref{numerical}. At the end, Section \ref{sec:conclusion} includes some concluding statements. 

\vspace{2mm}
\noindent{\textbf{Notation.}} Throughout the paper,  we adopt the following notations. Let $ \mathbb{C}^n$ and  $\R^{m\times n}$ denote the collection of all  $n$-dimensional complex vectors and  $m\times n$ real matrices, respectively. The notations $\|x\|_2,$ $\|A\|_F,$ $\sigma(A),$ $ {\vartheta(A)}$ and $tr(A)$ denote the $2$-norm,  {Frobenius} norm, spectrum,  {spectral radius} and trace of the matrix $A\in \R^{m\times n},$ respectively. For any vector $x\in \mathbb{C}^n,$ ${x}^*$ represents its conjugate transpose. For any matrix $A\in \R^{m\times n}$ having real eigenvalues, $\lambda_{\max}(A)$ and $\lambda_{\min}(A)$  denote the maximum eigenvalue and  minimum eigenvalue of $A,$ respectively.   {For any $x\in \mathbb{C}^n,$ $\mathfrak{R}(x)$ and $\Im(x)$ represent the real and imaginary part of $x,$ respectively. We denote $A_H:=\frac{A+A^T}{2}$ and $D_H:=\frac{D+D^T}{2}$.} 
\section{Solvability conditions for the DSPP and properties of the matrix $\mathcal{B}$}\label{sec2}
In this section, we provide the solvability conditions on the block matrices $A,$ $B,$ $C$ and $D$ for the DSPP \eqref{SPP1} and a few important properties of the matrix $\B.$
 {The nonsymmetric coefficient matrix $\mathcal{B}$ possesses the following desirable properties, which are crucial in the theoretical analysis of the iterative method and preconditioners designed to solve DSPP \eqref{SPP1}.
\begin{proposition}
   Let $A\in \R^{n\times n}$ and $D\in \R^{l\times l}$ with $A_H$ and $D_H$ be positive semidefinite. If $B$ has full row rank, then
   \begin{itemize}
       \item[(i)] $\B$ is semipositive real: ${\bm u}^T\B{\bm u}\geq 0$ for all ${\bm u}\in \R^{n+l+m}.$
       \item[(ii)] $\B$ is positive semistable: the real part of all eigenvalues of $\B$ is nonnegative.
   \end{itemize} 
\end{proposition}}
{ \proof $(i)$ Let $\bm{u}=[x^T,\,y^T\,z^T]^T\in \R^{n+l+m}.$ Then $\bm{u}^T\B\bm{u}=x^TAx+y^TDy.$ Therefore, we have $\bm{u}^T\B\bm{u}+\bm{u}^T\B^T\bm{u}=2(x^TA_Hx+y^TD_Hy)$ and hence, $\bm{u}^T\B\bm{u}=x^TA_Hx+y^TD_Hy\geq 0,$ as $A_H$ and $D_H$ are positive semidefinite. Thus, $\B$ is semipositive real.

\noindent$(ii)$ Let $\lambda$ be an eigenvalue of $\mathcal{B}$ and ${\bm u}=[u^T,v^T,w^T]^T\in \R^{n+l+m}$ is the corresponding eigenvector. Then $\bm{u}^*\B\bm{u}=\lambda \|\bm{u}\|_2$ and $(\bm{u}^*\B\bm{u})^*=\Bar{\lambda} \|\bm{u}\|_2.$ Thus 
\begin{align*}
 \mathfrak{R}(\lambda)&=   \frac{\bm{u}^*(\B+\B^T)\bm{u}}{2\|\bm{u}\|_2}\\
 &= \frac{\mathfrak{R}(\bm{u})^T(\B+\B^T)\mathfrak{R}(\bm{u})^T+\mathfrak{I}(\bm{u})^T(\B+\B^T)\mathfrak{I}(\bm{u})^T}{2\|\bm{u}\|_2}.
\end{align*}
Then, using $(i),$ we have $\mathfrak{R}(\lambda)\geq0.$ $\blacksquare$
}
\begin{proposition}\label{pro1}
   Let $A\in \R^{n\times n}$ and $D\in \R^{l\times l}$ be  {nonsingular matrices with $A_H$ and $D_H$ positive definite}. If $B$ has full row rank, then the saddle point matrix  $\mathcal{B}$ is nonsingular.
\end{proposition}
\proof Let $B$ has full row rank,  and ${\bm u}=[x^T,y^T,z^T]^T\in \R^{n+l+m}$ be such that  $\mathcal{B}{\bm u}=0.$ Then, we have
\begin{align}\label{eq21}
\left\{ \begin{array}{rc}
  Ax+B^Tz=&\0,\\
    Dy+Cz=&\0,\\ 
    -Bx-C^Ty=&\0.
    \end{array}\right.
\end{align}
We first assert that $x=\0.$ Then from the first equation in \eqref{eq21}, we obtain $B^Tz=\0.$ Since, $B$ has full row rank, this implies $z=\0.$ Thus the second equation of \eqref{eq21} gives $y=\0$ as $D$ is  {nonsingular}. This implies ${\bm u}=\0.$ Next, we assert $z=\0.$ Then, from the first and second  equation of \eqref{eq21}, we find that  $x=\0$ and  $y=\0,$ respectively, as $A$ and $D$ are  {nonsingular}. Hence, ${\bm u}=\0.$ Now, we assume that $x\neq \0$ and $z\neq \0.$ Then multiplying by $x^T$ from the left side of the first equation of  \eqref{eq21}, we obtain
\begin{equation}\label{eq24}
    x^TAx+x^TB^Tz=0.
\end{equation}
Again, multiplying third equation of \eqref{eq21}  by $z^T$ from the left, we get 
\begin{equation}\label{eq25}
    -z^TBx-z^TC^Ty=0.
\end{equation}
  Substituting \eqref{eq24} and  $Dy=-Cz$ on \eqref{eq25}, we have $x^TAx+y^TD^Ty=0,$   {and therefore it must be $x^TAx=0$ and $y^TD^Ty=0$, as both the quantities are nonnegative. However, $x^TAx=x^TA_Hx=0$ and $y^TD^Ty=y^TD_Hy=0$, which implies $x=\0$ and $y=\0,$ since $A_H$ and $D_H$ are positive definite matrices}. Thus, ${\bm u}=\0,$ and hence, $\mathcal{B}$ is nonsingular.  
$\blacksquare$

\begin{proposition}
    Let $A\in \R^{n\times n}$ be  symmetric positive definite (SPD) matrix and $D\in \R^{l\times l}$ be symmetric positive semidefinite matrix. If   $B$ and $C$ are of full row rank, then $\mathcal{B}$ is nonsingular.
\end{proposition}
\proof The proof follows similarly to the proof of proposition \ref{pro1}.

\begin{proposition}\label{prop3}
    Let $A\in \R^{n\times n}$ and $D\in \R^{l\times l}$  {be such that $A_H$ and $D_H$ are positive definite matrices}. If   $B$ and $C$ are of full row rank, then the matrix $\mathcal{B}$ is positive stable, i.e., $\mathfrak{R}(\lambda)>0$ for all $\lambda \in \sigma(\B).$ 
\end{proposition}
 \proof Suppose $\lambda$ is an eigenvalue of $\mathcal{B}$ and ${\bm u}=[u^T,v^T,w^T]^T\in \R^{n+l+m}$ is the corresponding eigenvector. Then, we have $\mathcal{B}{\bm u}=\lambda {\bm u},$ which leads to following three linear system of equations:
 \begin{align}\label{eq26}
 \left\{ \begin{array}{rc}
 Au+B^Tw     &=\lambda u,  \\
   Dv+Cw   & =\lambda v,\\
   -Bu-C^Tv&=\lambda w.
 \end{array}\right.
 \end{align}
Premultiplying $\mathcal{B}{\bm u}=\lambda {\bm u}$ by ${\bm u}^{*},$ we get 
\begin{equation}\label{eq27}
    \lambda \|{\bm u}\|_2^2=u^{*}Au+v^{*}Dv+2\iu \Im(u^{*}B^Tw+v^{*}Cw).
\end{equation}
 {On the other hand, from $\bm{u}^*\B^T\bm{u}=\Bar{\lambda}\|u\|_2^2,$ we get 
\begin{equation}\label{r1:eq27}
    \bar{\lambda} \|{\bm u}\|_2^2=u^{*}A^Tu+v^{*}D^Tv-2\iu \Im(u^{*}B^Tw+v^{*}Cw).
\end{equation}
By adding \eqref{eq27} and \eqref{r1:eq27}, we obtain
\begin{align}\label{r1:eq2} \nonumber(\lambda+\bar{\lambda})\|\bm{u}\|_2^2&=u^{*}(A+A^T)u+v^{*}(D+D^T)v\\ \nonumber
&=\mathfrak{R}(u)^T(A+A^T)\mathfrak{R}(u)+\mathfrak{I}(u)^T(A+A^T)\mathfrak{I}(u)\\
&+\mathfrak{R}(v)^T(D+D^T)\mathfrak{R}(v)+\mathfrak{I}(v)^T(D+D^T)\mathfrak{I}(v)
\end{align} 
Therefore, from \eqref{r1:eq2}, we obtain
\begin{equation}\label{eq28}
    \mathfrak{R}(\lambda)=\frac{\mathfrak{R}(u)^TA_H\mathfrak{R}(u)+\mathfrak{I}(u)^T A_H\mathfrak{I}(u)+\mathfrak{R}(v)^TD_H\mathfrak{R}(v)+\mathfrak{I}(v)^TD_H\mathfrak{I}(v)}{\|{\bm{u}}\|^2}.
\end{equation}}
 {Since $A_H$ and $D_H$ are positive definite matrices, from} \eqref{eq28}, we get $\mathfrak{R}(\lambda)\geq 0$ and $\mathfrak{R}(\lambda)=0$ if and only if $u=\0$ and $v=\0.$

 Next, we show that the vectors $u$ and $v$ can not be zero simultaneously. First, assume that $u=\0.$ From the first equation in \eqref{eq26}, it follows that $B^T w = \0$, which leads to $w=\0,$ as $B$ has full row rank.  Substituting $u = \0$ and $w = \0$ into the third equation of \eqref{eq26} yields $v = \0$.   As a result, we obtain \(\bm{u} = \mathbf{0}\), which contradicts the assumption that \(\bm{u}\) is an eigenvector. Therefore, we conclude that $u \neq \0$, which in turn implies that $\mathfrak{R}(\lambda) > 0$. This completes the proof. $\blacksquare$

 {To ensure the properties stated in Propositions \ref{pro1}-\ref{prop3} are satisfied, throughout the remainder of the paper, we assume \( A \) and \( D \) are nonsingular matrices, where \( A_H \) and \( D_H \) are positive definite matrices.}
\section{Proposed generalized shift-splitting (GSS) iterative method and preconditioner }\label{sec:GSS}
This section proposes a generalized shift-splitting (GSS) iterative method to solve the DSPP. Let  $\alpha, \beta, \tau,$ $\om$ be  positive real numbers and $P\in \R^{n\times n},Q\in \R^{l\times l},$ and $R\in \R^{m\times m}$ be SPD matrices, then the saddle point matrix $\mathcal{B}$ admits the following matrix splitting:
\begin{align}\label{GSS:1}
\mathcal{B}=  (\Theta+{\bm \omega}\B)-(\Theta-(1-\om) \B)=: {\P}_{\GSS}-\mathcal{N}_{\GSS},
\end{align}
where \begin{equation}\label{GSS}
     {\P}_{\GSS}=\bmatrix{\alpha P+\om A&\0& \om B^T\\ \0&\beta Q+\om D& \om C\\-\om B&-\om C^T& \tau R},
\end{equation}
\begin{equation}\label{QSS}
    \mathcal{N}_{\GSS}=\bmatrix{\alpha P-(1-\om )A&\0&(1-\om ) B^T\\ \0&\beta Q-(1-\om ) D & (1-\om ) C\\-(1-\om ) B&-(1-\om ) C^T& \tau R},
\end{equation}
 and $$\Theta=\bmatrix{\alpha P&\0&\0\\ \0& \beta Q&\0\\ \0&\0& \tau R}.$$

 The special matrix splitting in equation \eqref{GSS:1} introduces a novel iteration scheme, known as the GSS iterative method, for solving the DSPP.
 \begin{method}(GSS iterative method). Given the initial guess vector ${\bm u}^{(0)}=\left[{x^{(0)}}^T,{y^{(0)}}^T,{z^{(0)}}^T\right]^T,$  positive real numbers $\alpha, \beta, \tau$ and $\om,$  SPD matrices $P\in \R^{n\times n},Q\in \R^{l\times l}$ and  $R\in \R^{m\times m},$  until the stopping criterion is satisfied,  compute
 \begin{equation}\label{eq:iteration}
     {\bm u}^{(k+1)}=\mathcal{G}{\bm u}^{(k)}+{\bm d}, ~~ ~~k=0,1,2,\ldots,
 \end{equation}
     where ${\bm u}^{(k)}=\left[{x^{(k)}}^T,{y^{(k)}}^T,{z^{(k)}}^T\right]^T\in \R^{n+l+m},$ $\mathcal{G}= {\P}^{-1}_{\GSS}\mathcal{N}_{\GSS}$ is the iteration matrix and ${\bm d}= {\P}_{\GSS}^{-1}{\bm b}\in \R^{n+m+l}. $
 \end{method}
 The matrix splitting in \eqref{GSS:1} induces a preconditioner, denoted as $ {\P}_{\GSS}$, which can be utilized to speed up the convergence rate of the Krylov subspace iterative methods, like GMRES. This preconditioner is referred to as the GSS preconditioner. 
 
 At each step of the GSS iterative method or GSS preconditioned GMRES (PGMRES) method, we are required to solve a system of  linear equations in the following form: 
 \begin{equation}
      {\P}_{\mathrm{\GSS}}{\bm z}=r,
 \end{equation}
 where ${\bm z}=[{\bm z}^T_1,{\bm z}^T_2,{\bm z}^T_3]^T\in \R^{n+l+m}$ and ${\bm r}=[{\bm r}_1^T,{\bm r}_2^T,{\bm r}_3^T]^T\in \R^{n+l+m}.$ However, $ {\P}_{\GSS}$ admits the following decomposition:
 \begin{align}\label{decom:GSS}
   \nonumber   {\P}_{\GSS}=&\bmatrix{I& \0& \0\\ \0&I &\0 \\-\om B(\alpha P+\om A)^{-1}& -\om C^T(\beta Q+\om D)^{-1}&I}\bmatrix{\alpha P+\om A&\0&\0\\ \0& \beta Q+\om D&\0\\ \0 &\0&\widehat{R}}\\
     &\hspace{4cm}.\bmatrix{I&\0 &\om(\alpha P+\om A)^{-1}B^T \\ \0&I&\om (\beta Q+\om D)^{-1}C\\ \0&\0&I},
 \end{align}
 where $\widehat{R}=\tau R+\om^2 B(\alpha P+\om A)^{-1}B^T+\om^2 C^T(\beta Q+\om D)^{-1}C.$
 In the following, we present the algorithmic implementations of the GSS preconditioner designed to accelerate the GMRES method.
 \begin{algorithm}
			\caption{Solving $ {\P}_{\GSS}{\bm z}={\bm r}$}
                \label{alog1}
			\begin{algorithmic}
				\State \textbf{Input:} The matrices $A\in \R^{n\times n},\,B\in \R^{m\times n},\,C\in \R^{l\times m},\,D\in \R^{l\times l},$ $\bm{r}\in \R^{n+l+m},$ positive parameters $\alpha, \beta, \tau, \om,$ SPD matrices $P\in \R^{m\times m},$ $Q\in \R^{l\times l}$ and $R\in \R^{m\times m}.$ 
                \State \textbf{Output:} Solution vector ${\bm z}=[{\bm z}_1^T,{\bm z}_2^T, {\bm z}_3^T]^T\in \R^{n+l+m}.$
                \State \textbf{Steps:}
                \State  $1:$  Solve $(\alpha P+\om A)t_1={\bm r}_1$ to find $t_1.$
                \State  $2:$  Solve $( \beta Q+\om D)t_2={\bm r}_2$ to find $t_2.$
               \State  $3:$  Solve $\widehat{R}{\bm z}_3={\bm r}_3+\om Bt_1+\om C^Tt_2$ to obtain ${\bm z}_3.$
               \State  $4:$  Solve $(\alpha P+\om A){\bm z}_1 = {\bm r}_1-B^T {\bm z}_3$ to find ${\bm z}_1.$
               \State  $5:$  Solve $(\beta Q+\om D){\bm z}_2 = {\bm r}_2-\om C {\bm z}_3$ to obtain ${\bm z}_2.$
			\end{algorithmic}
		\end{algorithm}
	\begin{remark}
	    Algorithm \ref{alog1} necessitates solving two linear subsystems having the coefficient matrix $(\alpha P + \omega A)$,  two subsystems having the coefficient matrix $(\beta Q+\om D),$  and one subsystem with the coefficient matrix $\widehat{R}.$  {We can use LU factorization to solve them efficiently. Moreover, when $A$ and $D$ are SPD matrices, we have the flexibility to employ exact solvers, such as Cholesky factorization, and inexact solvers, like the preconditioned conjugate gradient method, to solve them efficiently}. Nevertheless, to solve steps $1$ and $4$, only one Cholesky  {or LU} factorization of \( \alpha P + \om A \) and to solve steps $2$ and $5$, only one Cholesky  {or LU} factorization of \( \beta Q + \om D \) is needed to perform.
	\end{remark}
  { \begin{remark}\label{remark2}
From Algorithm \ref{alog1}, observe that the most tedious task to implement the GSS preconditioner is to solve the linear subsystem in step $3.$   To avoid the direct construction of the matrices $B(\alpha P+\om A)^{-1}B^T$ and $C^T(\beta Q+\om D)^{-1}C,$ we modify the decomposition in \eqref{decom:GSS} in the following way: 
\begin{align}\label{decom_GSS2}
  \nonumber    {\widetilde{\P}}_{\GSS}:=&\bmatrix{I& \0& \0\\ \0&I &\0 \\-\om B(\alpha P+\om A)^{-1}& -\om C^T(\beta Q+\om D)^{-1}&I}\bmatrix{\alpha P+\om A&\0&\0\\ \0& \beta Q+\om D&\0\\ \0 &\0& \widetilde{P}+\widetilde{Q}}\\
     &\hspace{4cm}.\bmatrix{I&\0 &\om(\alpha P+\om A)^{-1}B^T \\ \0&I&\om (\beta Q+\om D)^{-1}C\\ \0&\0&I},
 \end{align}
 where $\widetilde{P}$ and $\widetilde{Q}$ are (efficient and economical) approximations of the matrices $B(\alpha P+\om A)^{-1}B^T$ and $C^T(\beta Q+\om D)^{-1}C,$ respectively.  Then, the step $4$ in Algorithm \ref{alog1}  changes to the linear subsystem $(\widetilde{P}+\widetilde{Q})\bm{z}_3={\bm r}_3+\om Bt_1+\om C^Tt_2.$ With suitably chosen $\widetilde{P}$ and $\widetilde{Q}$ (see for Example $1$), this subsystem is much easier to implement than step $3$ of Algorithm \ref{alog1}. We denote this inexact version of the GSS preconditioner as $\widetilde{P}_{\text{GSS}}.$
  \end{remark}}
  \section{Convergence analysis of the GSS iterative method}\label{SEC:convergence}
  The purpose of this section is to investigate the convergence behavior of the proposed GSS iterative method. To achieve this, the following result plays a crucial role.
  \begin{lemma}\label{lm1}
     Let $A\in \R^{n\times n}$ and $D\in \R^{l\times l}$  {be such that $A_H$ and $D_H$ positive definite matrices}, $B\in \R^{m\times n}$ and $C\in \R^{l\times m}$  {be} full row matrices. Then the matrix $\Theta^{-1}\mathcal{B}$ is positive stable. 
   \end{lemma}
  
  \proof Since $\Theta$ is SPD, $\Theta^{-1}\mathcal{B}$ is similar to the matrix $\Theta^{-\frac{1}{2}}\mathcal{B}\Theta^{-\frac{1}{2}}.$ By computation, we find that the block structure of the matrix \(\Theta^{-\frac{1}{2}}\mathcal{B}\Theta^{-\frac{1}{2}}\) is identical to that of \(\mathcal{B}\). Therefore, using Proposition \ref{prop3}, it follows that $\Theta^{-\frac{1}{2}}\mathcal{B}\Theta^{-\frac{1}{2}},$ and hence, $\Theta^{-1}\mathcal{B}$ is positive stable. $\blacksquare$

As noted in \cite{YSAAD}, a stationary iterative method in the form \eqref{eq:iteration} converges if and only if the iteration matrix has the spectral radius strictly less than one. The following result discusses the convergence of the GSS iterative method \eqref{eq:iteration}.
  \begin{theorem}\label{th1}
      Assume that $A\in \R^{n\times n}$ and $D\in \R^{l\times l}$  {are nonsingular matrices with $A_H$ and $D_H$ positive definite}, $B\in \R^{m\times n}$ and $C\in \R^{l\times m}$ are full row matrices. Let $\alpha, \beta, \tau, \om>0,$  $P\in \R^{n\times n}, Q\in \R^{l\times l}$ and $R\in \R^{m\times m}$ be SPD matrices. Then, the GSS iterative method converges to the unique solution of the DSPP  \eqref{SPP1} if   $$\dm{\om \geq \max\left\{\frac{1}{2}-\frac{\lambda_{\min}(\widetilde{\Theta})}{\vartheta(\Theta^{-1}\mathcal{B})^2},0\right\}},$$ where $\widetilde{\Theta}=\dm{\frac{\Theta^{-1}\mathcal{B}+\mathcal{B}^T\Theta^{-1}}{2}}.$
      
  \end{theorem}
  \proof The iteration matrix of the GSS iterative method \eqref{eq:iteration} is
  \begin{align}
\mathcal{G}= {\P}_{\GSS}^{-1}\mathcal{N}_{\GSS}&=(\Theta+\om \mathcal{B})^{-1}(\Theta-(1-\om) \mathcal{B})\\
&=(I {+}\,\om \Theta^{-1}\mathcal{B})^{-1}(I-(1-\om)\Theta^{-1}\mathcal{B}).
  \end{align}
Let $\lambda$  be an eigenvalue of  $\mathcal{G}$. Then
  \begin{equation}\label{eq44}
      (I {+}\,\om \Theta^{-1}\mathcal{B})^{-1}\left(I-(1-\om)\Theta^{-1}\mathcal{B}\right){\bf x}=\lambda {\bf x},
  \end{equation}
  where ${\bf x}\in \R^{n+m+l}$ is the corresponding eigenvector.
  From \eqref{eq44}, we can write
  \begin{align}
    &(I-(1-\om)\Theta^{-1}\mathcal{B})   {\bf x}=\lambda(I {+}\,\om \Theta^{-1}\mathcal{B}) {\bf x}\\ \label{eq46}
    \implies &((1-\om)+\lambda\om)\Theta^{-1}\mathcal{B}{\bf x}=(1-\lambda){\bf x}.
   \end{align}
   Note that $\lambda\neq 1,$ otherwise  \eqref{eq46} reduces to $\Theta^{-1}\mathcal{B}{\bf x}=0,$ which implies that ${\bf x}=0.$
   On the other hand, $(1-\om)+\lambda\om\neq 0,$ otherwise $(1-\lambda ){\bf x}=0.$ This gives ${\bf x}=0,$ which contradicts the assumption that ${\bf x}$ is an eigenvector. Therefore, from \eqref{eq46} we get 
   \begin{equation}
      \Theta^{-1} \mathcal{B}{\bf x}=\frac{1-\lambda}{(1-\om)+\lambda\om}{\bf x}.
   \end{equation}
   Thus $\theta:=\frac{1-\lambda}{(1-\om)+\lambda\om}$ is an eigenvalue of $\Theta^{-1} \mathcal{B}.$ Further, we can write $$\lambda=\frac{1-(1-\om)\theta}{1+\om \theta}.$$

   \noindent Therefore, $|\lambda|<1$ if and only if 
   $ |1-(1-\om)\theta|< |1+\om \theta|,$
\noindent   i.e.,
   \begin{equation}\label{eq48}
 (1-(1-\om)\mathfrak{R}(\theta))^2+(1-\om)^2\Im(\theta)^2<(1+\om \mathfrak{R}(\theta))^2+\om^2\Im(\theta).
   \end{equation}
   \noindent Consequently, it follows from \eqref{eq48} that the iterative method \eqref{eq:iteration} is convergent if
   \begin{equation}\label{eq49}
       2 \mathfrak{R}(\theta)+(2\om-1)|\theta|^2>0.
   \end{equation}
   By Lemma \ref{lm1}, we have $\mathfrak{R}(\theta)>0,$ and this implies $\frac{\mathfrak{R}(\theta)}{|\theta|^2}>0.$ From \eqref{eq49}, we get $\om>\frac{1}{2}-\frac{\mathfrak{R}(\theta)}{|\theta|^2}.$
   
 \noindent Next, assume that $\bm{w}$ is the eigenvector corresponding to the eigenvalue $\theta.$ Then $\Theta^{-1}\mathcal{B}\bm{w}=\theta \bm{w}.$ Multiplying by $\bm{w}^*$ from the left side, we have $\bm{w}^*\Theta^{-1}\mathcal{B} \bm{w} =\theta \bm{w}^* \bm{w},$ and taking conjugate transpose  gives $\bm{w}^*\mathcal{B}^T\Theta^{-1}\bm{w}=\Bar{\theta} \bm{w}^*\bm{w}.$ Hence, $$\mathfrak{R}(\theta)=\frac{\bm{w}^*(\Theta^{-1}\mathcal{B}+\mathcal{B}^T\Theta^{-1})\bm{w}}{2\bm{w}^*\bm{w}}\geq \lambda_{\min}(\widetilde{\Theta}).$$
   Again $|\theta|\leq\vartheta(\Theta)$ gives $\dm{\frac{1}{2}-\frac{\mathfrak{R}(\theta)}{|\theta|^2}\leq \frac{1}{2}-\frac{\lambda_{\min}(\widetilde{\Theta})}{\vartheta(\Theta)^2}}.$
   Since $\om>0,$ the GSS iterative method is convergent if $\dm{\om >\max\left\{\frac{1}{2}-\frac{\lambda_{\min}(\widetilde{\Theta})}{\vartheta(\Theta)^2},0\right\}}.$
  $\blacksquare $
  \begin{remark}\label{rem1}
      Note that if $\om \geq \frac{1}{2},$ then the condition \eqref{eq48} holds. This shows that the GSS iterative method \eqref{eq:iteration} is convergent for any initial guess vector if $\om \geq \frac{1}{2}.$
  \end{remark}
Notably, solving the DSPP \eqref{SPP1} is same as to solve the preconditioned linear system \(\P_{\text{GSS}}^{-1}\B {\bf u} = {\P}_{\text{GSS}}^{-1} {\bm b}\). Hence, as an immediate consequence of Theorem \ref{th1} and Remark \ref{rem1}, we have the following results regarding the clustering properties of the spectrum of the preconditioned matrix  $ {\P}_{\text{GSS}}^{-1}\B$.
  \begin{theorem}
    Assume that $A\in \R^{n\times n}$ and $D\in \R^{p\times p}$ {  are nonsingular matrices with $A_H$ and $D_H$ positive definite},  $B\in \R^{m\times n}$ and $C\in \R^{p\times m}$ are full row matrices. Let $ \P_{\mathrm{GSS}}$ be defined as in \eqref{GSS} and $\lambda$ be an eigenvalue of the preconditioned matrix $ \P_{\mathrm{GSS}}^{-1}\B$. If $\om\geq 1/2,$ then $\lambda$ satisfies the following:
    \begin{eqnarray*}
        |\lambda-1|<1,
    \end{eqnarray*}
    i.e., all eigenvalues of the preconditioned matrix $ \P_{\mathrm{GSS}}^{-1}\B$ are entirely contained in a circle centered at $(1, 0)$ with radius strictly less one.
\end{theorem}
\proof The proof follows immediately from the identity $ \P_{\text{GSS}}^{-1}\B= \P_{\text{GSS}}^{-1}( \P_{\text{GSS}}-\mathcal{N}_{\text{GSS}})=I-\mathcal{G}.$ $\blacksquare$

  \section{ Two relaxed variants of GSS preconditioner}\label{sec:RGSS}
To enhance efficiency, this section introduces two relaxed variants of the GSS preconditioner. By removing $\alpha P$ from $(1,1)$ block  and $\beta Q$ from  $(2,2)$ block of $ \P_{\GSS},$ we obtain the following two relaxed GSS (RGSS) preconditioners:
  \begin{equation}
       \P_{\text{RGSS-I}}=\bmatrix{\om A &\0& \om B^T\\ \0& \beta Q+\om D& \om C\\-\om B& -\om C^T& \tau R}
  \end{equation}
  and 
   \begin{equation}
       \P_{\text{RGSS-II}}=\bmatrix{\om A &\0& \om B^T\\ \0&\om D& \om C\\-\om B& -\om C^T& \tau R}.
  \end{equation}

  In the implementation, at each step of RGSS-I and RGSS-II preconditioners in conjunction with GMRES, we are required to solve the linear systems of the following forms: 
  \begin{equation}
       \P_{\text{RGSS-I}}{\bm w}={\bm r}~~ \text{or}~~   \P_{\text{RGSS-II}}{\bm w}={\bm r},
  \end{equation}
 respectively.
 
  Set ${\mathcal{R}}_1=\tau R+\om BA^{-1}B^T+\om^2 C^T(\beta Q+\om D)^{-1}C$ and $\mathcal{R}_2=\tau R+\om BA^{-1}B^T+\om C^T D^{-1}C,$ then 
\begin{align}\label{decom:GSSI}
      \P_{\text{RGSS-I}}=&\bmatrix{I& \0& \0\\ \0&I &\0 \\- BA^{-1}& - \om C^T (\beta Q+\om D)^{-1} & I}\bmatrix{ \om A&\0&\0\\ \0& \beta Q+\om D&\0\\ \0 &\0&\mathcal{R}_1}\bmatrix{I&\0 & A^{-1}B^T \\ \0 &I & \om (\beta Q+\om D)^{-1} C \\ \0&\0 & I},
 \end{align}
\begin{align}\label{decom:GSSII}
      \P_{\text{RGSS-II}}=&\bmatrix{I& \0& \0\\ \0&I &\0 \\- BA^{-1}& - C^T D^{-1} & I} \bmatrix{ \om A&\0&\0 \\ \0&  \om D & \0 \\ \0 &\0&\mathcal{R}_2}\bmatrix{I&\0 & A^{-1}B^T \\ \0 &I & D^{-1} C \\ \0&\0 & I},
 \end{align}
  By applying a similar methodology as in Algorithm \ref{alog1}, we derive Algorithms \ref{alog2} and \ref{alog3} for implementing the RGSS-I and RGSS-II preconditioners, respectively.
 
   \begin{algorithm}
			\caption{Solving $ \P_{\text{RGSS-I}}{\bm w}={\bm r}$}
                \label{alog2}
			\begin{algorithmic}
				\State \textbf{Input:}  The matrices $A\in \R^{n\times n},\,B\in \R^{m\times n},\,C\in \R^{l\times m},\,D\in \R^{l\times l},$ $\bm{r}\in \R^{n+l+m},$ positive parameters $ \beta, \tau, \om,$ SPD matrices  $Q\in \R^{l\times l}$ and $R\in \R^{m\times m}.$  
                \State \textbf{Output:} Solution vector ${\bm w}=[\w_1^T, \w_2^T, \w_3^T]^T\in \R^{n+l+m}.$
                \State \textbf{Steps:}
                \State  $1:$  Solve $A t_1=\bm{r}_1/\om $ to find $t_1.$
                \State  $2:$  Solve $( \beta Q+\om D)t_2=\bm{r}_2$ to find $t_2.$
               \State  $3:$  Solve $\mathcal{R}_1 \w_3=\bm{r}_3+\om Bt_1+\om C^Tt_2$ to obtain $\w_3.$
               \State  $4:$  Solve $ A \w_1=\frac{1}{\om}(\bm{r}_1-B^T\w_3)$ to find $\w_1.$
               \State  $5:$  Solve $(\beta Q+\om D)\w_2=\bm{r}_2-\om C \w_3$ to obtain $\w_2.$
			\end{algorithmic}
		\end{algorithm}
 \begin{algorithm}
			\caption{Solving $ \P_{\text{RGSS-II}}{\bm w}={\bm r}$}
                \label{alog3}
			\begin{algorithmic}
				\State \textbf{Input:} The matrices $A\in \R^{n\times n},\,B\in \R^{m\times n},\,C\in \R^{l\times m},\,D\in \R^{l\times l},$ $\bm{r}\in \R^{n+l+m},$ positive parameters $ \tau, \om$ and the SPD matrix  $R\in \R^{m\times m}.$ 
                \State \textbf{Output:} Solution vector ${\bm w}=[\w_1^T,\w_2^T, \w_3^T]^T\in \R^{n+l+m}.$
                \State \textbf{Steps:}
                \State  $1:$  Solve $ A t_1=\bm{r}_1/\om$ to find $t_1.$
                \State  $2:$  Solve $ D t_2=\bm{r}_2/\om$ to find $t_2.$
               \State  $3:$  Solve $\mathcal{R}_2\w_3=\bm{r}_3+\om B t_1+\om C^Tt_2$ to obtain $\w_3.$
               \State  $4:$  Solve $ A \w_1=\frac{1}{\om}( \bm{r}_1-B^T\w_3)$ to find $\w_1.$
               \State  $5:$  Solve $ D \w_2=\frac{1}{\om} (\bm{r}_2-\om C \w_3)$ to obtain $\w_2.$
			\end{algorithmic}
		\end{algorithm}
{ \begin{remark}
A key challenge in implementing the RGSS-I and RGSS-II preconditioners lies in solving the linear subsystems associated with the coefficient matrices \(\mathcal{R}_1\) and \(\mathcal{R}_2\), respectively. As noted in Remark \ref{remark2},  to avoid the direct computation $BA^{-1}B^T$ and $C^T(\beta Q+\om D)^{-1}C$ in Algorithm \ref{alog2}, and  $BA^{-1}B^T$ and $C^T D^{-1}C$ in Algorithm \ref{alog3}, we can use approximate versions of these terms. Following a similar technique as in decomposition \eqref{decom_GSS2}, let $\widetilde{P}_1,$ $\widetilde{Q}_1$ and $\widetilde{Q}_2$ be efficient and economical approximations of the matrices  $B A^{-1}B^T,$ $C^T(\beta Q+\om D)^{-1}C,$ and $C^T D^{-1}C,$ respectively.  With these approximations, step $4$ in Algorithms  \ref{alog2} and \ref{alog3}   transforms into solving the linear subsystems $(\widetilde{P}_1+\widetilde{Q}_1)\bm{z}_3={\bm r}_3+\om Bt_1+\om C^Tt_2$ and   $(\widetilde{P}_1+\widetilde{Q}_2)\bm{z}_3={\bm r}_3+\om Bt_1+\om C^Tt_2,$ respectively.  By appropriately selecting $\widetilde{P}_1$ $\widetilde{Q}_1,$ and $\widetilde{Q}_2$ (see for Example $1$), these subsystems become significantly easier to implement compared to step $3$  of Algorithms \ref{alog2} and \ref{alog3}. The resulting inexact preconditioners are denoted by $\widetilde{\P}_{\text{RGSS-I}}$ and $\widetilde{\P}_{\text{RGSS-II}},$ respectively. 
\end{remark}}
 Next, we investigate the spectral properties of the preconditioned matrices $  \P_{\text{RGSS-I}}^{-1}\mathcal{B}$ and $  \P_{\text{RGSS-II}}^{-1}\mathcal{B},$ by considering $A$ and $D$ are SPD. 
\begin{theorem}\label{theorem51}
Assume that $A$ and $D$ are SPD matrices, $B$ and $C$ have full row rank, and let $Q$ and $R$ be SPD matrices. Then, the preconditioned matrix  $ \P_{\mathrm{RGSS}{\text-}\mathrm{I}}^{-1} \B$ has   $\dm{\frac{1}{\om}}$ as the eigenvalue with multiplicity $n.$ Further, let $\mu$ be an eigenvalue among the remaining $m+l$ eigenvalues with the corresponding eigenvector $[u^T, v^T]^T$ such that $\|\sqrt{\beta}Q^{1/2}u^T, \sqrt{\tau}R^{1/2}v^T]^T\|_2=1$. Then
\begin{enumerate}
    \item[(1)] if $\Im(\frac{\mu}{1-\om \mu}) \neq 0,$ we have $\beta\|Q^{1/2}u\|^2_2=\frac{1}{2}=\tau\|R^{1/2}v\|^2_2$ and 
    \begin{eqnarray*}
        \mathfrak{R}\left(\frac{\mu}{1-\om \mu}\right) = \frac{1}{2}\left(\frac{u^*Du}{\beta u^*Qu}+ \frac{v^*Sv}{\tau v^*Rv} \right).
    \end{eqnarray*}
    Thus, it satisfies the following bounds:
    \begin{align}\label{bound1:th51}
    \left\{\begin{array}{c}
        \frac{1}{2} \left(\frac{\lambda_{\min}(Q^{-1}D)}{\beta}+ \frac{\lambda_{\min}(R^{-1}S)}{\tau} \right) \leq \mathfrak{R}\left(\frac{\mu}{1-\om \mu}\right) \leq \frac{1}{2} \left(\frac{\lambda_{\max}(Q^{-1}D)}{\beta}+ \frac{ \lambda_{\max}(R^{-1}S) }{\tau}\right),\\
          |\Im\left(\frac{\mu}{1-\om \mu}\right)|\leq \sigma_{\max}\left(\frac{1}{\sqrt{\beta \tau}}R^{-\frac{1}{2}} C Q^{-\frac{1}{2}}\right).
    \end{array}\right.
    \end{align}
    \item[(2)] If $\Im(\frac{\mu}{1-\om \mu}) =0,$ we have $$\frac{\mu}{1-\om \mu}=\frac{u^* D u+ v^* S v}{\beta u^* Q u+\tau    v^* R v}$$
    and it holds that: 
     \begin{align}\label{bound2:th51}
    \begin{array}{c}
        2 \min\left\{\frac{\lambda_{\min}(Q^{-1}D)}{\beta}, 
 \frac{\lambda_{\min}(R^{-1}S)}{\tau} \right\}\leq \frac{\mu}{1-\om \mu} \leq \max\left\{\frac{\lambda_{\max}(Q^{-1}D)}{\beta}, 
 \frac{\lambda_{\max}(R^{-1}S)}{\tau}\right\}.
    \end{array}
    \end{align}
\end{enumerate}
\end{theorem}
\proof Let $S=BA^{-1}B^T.$ Then the  saddle point matrix $\B$  and the preconditioner $ \P_{\text{RGSS-I}}$ admits the following decompositions:
\begin{eqnarray}\label{dec1}
    \B=\mathbf{L}\mathbf{D}\mathbf{U} \quad \text{and} \quad  \P_{\text{RGSS-I}} = \mathbf{L} \widetilde{\mathbf{D}} \mathbf{U},
\end{eqnarray}
 where 
 \begin{align*}
     \mathbf{L}&=\bmatrix{I & \0 & \0 \\ \0 & I & \0\\ -BA^{-1} & \0 & I},~~\mathbf{D}=\bmatrix{A & \0 & \0 \\ \0 & D & C\\ \0 & -C^T & S},~~\widetilde{\mathbf{D}}=\bmatrix{\om A & \0 & \0 \\ \0 & \om D+ \beta Q & \om C\\ \0 & -\om C^T & S + \tau R}\\
    &\text{and}~~\mathbf{U} = \bmatrix{I & \0 & A^{-1}B^T \\ \0 & I & \0\\ \0 & \0 & I}.
 \end{align*} 
Using  decompositions in \eqref{dec1}, we obtain 
\begin{eqnarray}\label{thm51:2}
     \P_{\text{RGSS-I}}^{-1} \B =  \mathbf{U}^{-1} \widetilde{\mathbf{D}}^{-1} \mathbf{D} \mathbf{U}. 
\end{eqnarray}
Therefore, $ \P_{\text{RGSS-I}}^{-1} \B$ is similar to $\widetilde{\mathbf{D}}^{-1} \mathbf{D},$ which is given by
\begin{eqnarray}
   \widetilde{\mathbf{D}}^{-1} \mathbf{D} =  \bmatrix{\om^{-1}I & \0 \\ \0 & \mathbf{M}^{-1}{\mathbf{A}}},
\end{eqnarray}
 where $\mathbf{M}=\bmatrix{\om D+\beta Q & \om C \\ -\om C^T & \om S +\tau R} $ and $\mathbf{A}=\bmatrix{D & C \\ -C^T & S}.$ Hence,  $ \P_{\text{RGSS-I}}^{-1} \B$ has an eigenvalue $\dm{\frac{1}{\om }}$ with multiplicity at least $n$ and while the remaining eigenvalues are those of the preconditioned matrix $\mathbf{M}^{-1}\mathbf{A}.$ Consider $\mathbf{I}= \bmatrix{\0 & I \\ -I & \0 }\in\R^{l+m },$ then $\mathbf{M}^{-1}\mathbf{A}=\mathbf{I} \left( \mathbf{I}^{-1} \mathbf{M} \mathbf{I} \right)^{-1} \left(\mathbf{I}^{-1}\mathbf{A}\mathbf{I}\right) \mathbf{I}^{-1}.$ Consequently, $\mathbf{M}^{-1}\mathbf{A}$ is similar to $\widetilde{\mathbf{M}}^{-1} \widetilde{\mathbf{A}},$ where 
 \begin{eqnarray}
   \widetilde{\mathbf{M}}:=   \mathbf{I}^{-1} \mathbf{M} \mathbf{I} = \bmatrix{\om D+\beta Q & \om C^T \\ -\om C & \om S + \tau R} \quad \text{and} \quad \widetilde{\mathbf{A}}:=   \mathbf{I}^{-1} \mathbf{A} \mathbf{I} = \bmatrix{ D &  C^T \\ - C &  S}.
 \end{eqnarray}
  Let $\mu$ be an eigenvalue of the matrix $\widetilde{\mathbf{M}}^{-1} \widetilde{\mathbf{A}}$ with the corresponding eigenvalue $[u^T, v^T]^T.$ Assert that $\mu = \frac{1}{\om},$ then
 \begin{eqnarray}
     \bmatrix{\beta Q & \0 \\ \0 & \tau R}\bmatrix{u\\ v}= \0,
 \end{eqnarray}
 which gives $[u^T, v^T]^T=  \0.$ This contradicts to the  assumption that $[u^T, v^T]^T$ is an eigenvector, and hence $\mu \neq \frac{1}{\om}.$ Since, $\mu$ is an eigenvalue of $\widetilde{\mathbf{M}}^{-1} \widetilde{\mathbf{A}}$, we have 
 \begin{align*}
     &\mu \bmatrix{\om D+\beta Q & \om C^T \\ -\om C & \om S + \tau R} \bmatrix{u\\ v}= \bmatrix{D & C^T \\ -C & S}\bmatrix{u\\ v}\\
     &\Longrightarrow \frac{\mu}{1-\om \mu}\bmatrix{\beta  Q & \0 \\ \0 &\tau  R}\bmatrix{u\\ v}= \bmatrix{D & C^T \\ -C & S} \bmatrix{u\\ v}\\
     &\Longrightarrow \frac{\mu}{1-\om \mu} \bmatrix{\sqrt{\beta}Q^{\frac{1}{2}}u   \\ \sqrt{\tau} R^{\frac{1}{2}} v}= \bmatrix{\frac{1}{\sqrt{\beta}}Q^{-\frac{1}{2}} & \0 \\ \0 &\frac{1}{\sqrt{\beta}} R^{-\frac{1}{2}}}  \bmatrix{D & C^T \\ -C & S} \bmatrix{\frac{1}{\sqrt{\beta}}Q^{-\frac{1}{2}} & \0 \\ \0 &\frac{1}{\sqrt{\tau}} R^{-\frac{1}{2}}} \bmatrix{\sqrt{\beta}Q^{\frac{1}{2}} u \\ \sqrt{\tau}R^{\frac{1}{2}} v}\\
     &\Longrightarrow \frac{\mu}{1-\om \mu} \bmatrix{\widetilde{u}\\ \widetilde{v}} = \bmatrix{ \widetilde{D} &  \widetilde{C}^T \\ -\widetilde{C} &  \widetilde{S} } \bmatrix{ \widetilde{u}\\ \widetilde{v} },
 \end{align*}
 where $\widetilde{D}= \frac{1}{\beta}Q^{-\frac{1}{2}} D Q^{-\frac{1}{2}},$ $\widetilde{C}=\frac{1}{\sqrt{\beta \tau}}R^{-\frac{1}{2}} C Q^{-\frac{1}{2}},$ $\widetilde{S}= \frac{1}{\tau
 }R^{-\frac{1}{2}} S R^{-\frac{1}{2}},$ $\widetilde{u}= \sqrt{\beta} Q^{\frac{1}{2}} u$ and $\widetilde{v} =\sqrt{\tau} R^{\frac{1}{2}} v.$ Since, $\widetilde{D}$ and $\widetilde{S}$ are SPD matrices and $\widetilde{C}$ has full row rank, according to Proposition $2.12$ in \cite{Benzi2006}, we have the desired bounds of \eqref{bound1:th51} and \eqref{bound2:th51}. 
 $\blacksquare$

\begin{remark}
The asymptotic behavior of the eigenvalue \(\mu\) is analyzed according to Theorem \ref{theorem51} when the iteration parameters  $\beta$ and $\tau$ approach zero from the positive side. This analysis is conducted in the following two cases:
\begin{itemize}
    \item Let $\theta=\frac{\mu}{1-\om \mu}.$ When $\Im(\theta)\neq 0, $
 we have $\mathfrak{R}(\theta)\rightarrow +\infty$ and $|\Im(\theta)|\rightarrow +\infty $ as $\beta, \tau \rightarrow 0_{+}.$ Then, it holds that \begin{eqnarray*}
     \mu= \left[\frac{1}{\om}-\frac{1+\om \mathfrak{R}(\theta)}{(1+\om \mathfrak{R}(\theta))^2+\om^2 \Im(\theta)^2}\right]+\iu \frac{\Im(\theta)}{(1+\om \mathfrak{R}(\theta))^2+\om^2 \Im(\theta)^2}\rightarrow \frac{1}{\om} ~\text{as} ~ \beta, \tau \rightarrow 0_{+}. 
 \end{eqnarray*}
 \item When $\Im(\theta)= 0,$ we have $\theta \rightarrow + \infty$ as $\beta, \tau \rightarrow 0_{+}.$ Then, it holds that
 \begin{eqnarray}
     \mu=\frac{\theta}{1+ \om \theta} \rightarrow \frac{1}{\om}, ~\text{as}~   \beta, \tau \rightarrow 0_{+}.
 \end{eqnarray}
 \end{itemize}
\end{remark}

  The following theorem establishes the spectral properties for the RGSS-II preconditioned matrix  $  \P_{\mathrm{RGSS}\text{-}\mathrm{II}}^{-1}\mathcal{B}$.
 \begin{theorem}\label{th52}
     Assume that $A$ and $D$ are SPD matrices, $B$ and $C$ have full row rank, and let $R$ be an SPD matrix. Then the preconditioned matrix $  \P_{\mathrm{RGSS}\text{-}\mathrm{II}}^{-1}\mathcal{B}$ has the eigenvalue $\dm{\frac{1}{\om}}$ with multiplicity $n+l.$ The remaining eigenvalues satisfy the generalized eigenvalue problem $(BA^{-1}B^T+C^TD^{-1}C+ \tau R){\bf x}=\lambda \mathcal{R}_2{\bf x},$ where $\mathcal{R}_2=\tau R+\om BA^{-1}B^T+\om C^T D^{-1}C.$
  \end{theorem}
\proof From \eqref{decom:GSSII}, we have 
\begin{align}\label{eq1:th52}
  \nonumber   \P_{\text{RGSS-II}}^{-1}\B&= \bmatrix{\om^{-1}A & \0 & -A^{-1}B^T\mathcal{R}_2^{-1} \\ \0 & \om^{-1}D & -D^{-1}C\mathcal{R}_2^{-1}\\ \0 & \0 &\mathcal{R}_2} \bmatrix{A& \0& B^T \\ \0 &D & C\\ \0 & \0 & \mathbf{X}}\\
    &=\bmatrix{\om^{-1}I & \0 & \om^{-1}AB^T-A^{-1}B^T \mathcal{R}_2^{-1} \mathbf{X}\\ \0 & \om^{-1}I & \om^{-1}DC - D^{-1} C \mathcal{R}_2^{-1} \mathbf{X}\\ \0 & \0 &  \mathcal{R}_2^{-1} \mathbf{X} },
\end{align}
where $\mathbf{X}=BA^{-1}B^T+C^TD^{-1}C+ \tau R.$ Thus, $ \P_{\text{RGSS-II}}^{-1}\B$ has the eigenvalue $\dm{\lambda=\frac{1}{\om}}$ with multiplicity at least $n+l,$ and the rest of  eigenvalues satisfy the generalized eigenvalue problem 
\begin{equation*}
   (BA^{-1}B^T+C^TD^{-1}C+\tau R) {\bf x} =\lambda  \mathcal{R}_2 {\bf x}.
\end{equation*}
Hence, the proof is completed. $\blacksquare$

\begin{corollary}
    Suppose that assumption on Theorem \ref{th52} holds. Then, the eigenvalues of $  \P_{\mathrm{RGSS}\text{-}\mathrm{II}}^{-1}\mathcal{B}$  {satisfy} 
    \begin{equation}
        \lambda \in [\Lambda_{\min}, \Lambda_{\max}],
    \end{equation}
    where $$\Lambda_{\min}=\frac{ \eta_{\min}+ \xi_{\min} +\tau}{ \om \eta_{\max}+ \om \xi_{\max} +\tau},~ \Lambda_{\max}=\frac{ \eta_{\max}+ \xi_{\max} +\tau}{ \om \eta_{\min}+ \om \xi_{\min} +\tau},$$
    $\eta_{\min}=\lambda_{\min}(R^{-1}BA^{-1}B^T),$ $\eta_{\max}=\lambda_{\max}(R^{-1}BA^{-1}B^T),$ $\xi_{\min}= \lambda_{\min}(R^{-1}C^TD^{-1}C)$ and $\xi_{\max}=  \lambda_{\max}(R^{-1}C^TD^{-1}C).$ 
\end{corollary}  
\proof Premultiplying by ${\bf x}^T$ of the generalized eigenvalue problem $(BA^{-1}B^T+C^TD^{-1}C+\tau R){\bf x}=\lambda \mathcal{R}_2{\bf x},$ we obtain
\begin{align}\label{eq:coro51}
  \nonumber  \lambda &=\frac{{\bf x}^T(BA^{-1}B^T + C^TD^{-1}C + \tau R){\bf x}}{{\bf x}^T(\om BA^{-1}B^T+\om C^T D^{-1}C + \tau R){\bf x}}\\
    &=\frac{{(R^{1/2}\bf x)}^TR^{-1/2}(BA^{-1}B^T + C^TD^{-1}C + \tau)R^{-1/2}{(R^{1/2}\bf x)}}{{(R^{1/2}\bf x)}^TR^{-1/2}(\om BA^{-1}B^T +\om C^TD^{-1}C + \tau)R^{-1/2}{(R^{1/2}\bf x)}}.
\end{align}
Since  $R^{-1/2}BA^{-1}B^TR^{-1/2}$ is similar to $R^{-1}BA^{-1}B^T$ and $R^{-1/2}C^TA^{-1}CR^{-1/2}$ is similar to $R^{-1}C^TA^{-1}C,$ and $\lambda_{\min}(X)\leq \theta\leq \lambda_{\max}(X)$ for any $\theta \in \sigma(X),$ where $X$ is SPD, the proof follows from \eqref{eq:coro51}.
$\blacksquare$

  Next, we will examine the properties of the minimal polynomial of the preconditioned matrix  $  \P_{\text{RGSS-II}}^{-1}\mathcal{B},$  which determines the dimension of the Krylov subspace.
  \begin{theorem}\label{th53}
      Assume that $A\in \R^{n\times n}$ and $D\in \R^{l\times l}$ are nonsingular matrices with $A_H$ and $D_H$ are SPD matrices,  $B \in\R^{m\times n}$ and $C\in \R^{l\times m}$ are full row rank matrices. Then, the degree of the minimal polynomial of the preconditioned matrix $  \P_{\mathrm{RGSS-II}}^{-1}\mathcal{B}$ is at most $m+1.$ Therefore, the dimension of the Krylov subspace $\mathcal{K}( \P_{\mathrm{RGSS-II}}^{-1}\mathcal{B}, {\bm b})$ is at most $m+1.$
  \end{theorem}
\proof From \eqref{eq1:th52}, we obtain  
\begin{eqnarray}
     \P_{\text{RGSS-II}}^{-1}\mathcal{B}=\bmatrix{\om^{-1}I & \0& \Sigma_1\\ \0&\om^{-1}I& \Sigma_2\\ \0&\0&\Sigma_3},
\end{eqnarray}
where $\Sigma_1=\om^{-1}AB^T-A^{-1}B^T \mathcal{R}_2^{-1} \mathbf{X},$ $\Sigma_2=\om^{-1}DC - D^{-1} C \mathcal{R}_2^{-1} \mathbf{X}, $ $\Sigma_3= \mathcal{R}_2^{-1} \mathbf{X},$ and $\mathbf{X}=BA^{-1}B^T+C^TD^{-1}C+R.$
Let $\mu_i,$ $i=1,2,\ldots,m$ be the eigenvalues of  of $\Sigma_3.$ Then they  are also  eigenvalues of the preconditioned matrix $  \P_{\text{RGSS-II}}^{-1}\mathcal{B}.$ Then the characteristic polynomial of $  \P_{\text{RGSS-II}}^{-1}\mathcal{B}$ is given by
\begin{equation*}
    f(\lambda)=\left(\lambda-\frac{1}{\om}\right)^{n+l}\prod_{i=1}^m(\lambda-\mu_i).
\end{equation*}
Consider the polynomial ${ g(\lambda)=(\lambda-\frac{1}{\om})\prod\limits_{i=1}^{m}(\lambda-\mu_i)}.$ Then 
\vspace{-2mm}
\begin{align}
  \nonumber  g( \P_{\text{GSS-II}}^{-1}\mathcal{B})&=( \P_{\text{GSS-II}}^{-1}\mathcal{B}-\frac{1}{\om}I)\prod\limits_{i=1}^{m}( \P_{\text{GSS-II}}^{-1}\mathcal{B}-\mu_iI)\\ \label{eq:Th52}
            &=\bmatrix{\0&\0&\Sigma_1\\ \0&\0&\Sigma_2\\ \0&\0&\Sigma_3-\om^{-1}I}\bmatrix{{\prod\limits_{i=1}^{n}(\om^{-1}I-\mu_iI)}&\0&\Sigma_1\\ \0&{\prod\limits_{i=1}^{l}(\om^{-1}I-\mu_iI)}&\Sigma_2\\ \0&\0&{\prod\limits_{i=1}^{m}(\Sigma_3-\mu_iI)}}.
\end{align}
Given that $\Sigma_3$ has the eigenvalues $\mu_i,$ $i=1,2,\ldots,m,$   we obtain ${\prod\limits_{i=1}^{m}(\Sigma_3-\mu_iI)}=\0.$ Therefore, from \eqref{eq:Th52}, we obtain $g( \P_{\text{RGSS-II}}^{-1}\mathcal{B})=\0.$ Hence, by Cayley-Hamilton theorem, the degree of the minimal polynomial of $ \P_{\text{RGSS-II}}^{-1}\mathcal{B}$ is at most $m+1.$

\noindent As mentioned in \cite{YSAAD}, the degree of the minimal polynomial of a matrix and the dimension of the associated Krylov subspace are equal. Therefore, the  Krylov subspace $\mathcal{K}( \P_{\text{RGSS-II}}^{-1}\mathcal{B}, {\bm b})$ has dimension  at most $m+1.$ Hence, the proof follows.  $\blacksquare$

Based on the property outlined in Theorem \ref{th53}, Krylov subspace methods like PGMRES require maximum \(m+1\) iterations to solve the system \eqref{SPP1}.
\section{Discussion on the selection of the parameters}\label{sec:parameter}
It is worth noting that proposed GSS, RGSS-I and RGSS-II preconditioners involve the parameters $\alpha,$ $\beta,$ $\tau$ and $\om.$ As a general criterion for a preconditioner to perform efficiently, it should be as close as possible to the coefficient matrix of the system \cite{Benzi2005}, we find the optimal choices for the proposed GSS preconditioner by minimizing $\|\mathcal{N}_{\text{GSS}}\|_F=\| \P_{\text{GSS}}-\B\|_F.$ First, we define the function $\varphi$ by $$\varphi(\alpha, \beta, \tau, \om)= \|\mathcal{N}_{\text{GSS}}\|_F^2 = tr(\mathcal{N}_{\text{GSS}}^T\mathcal{N}_{\text{GSS}})>0.$$
After some easy calculations, we obtain
\begin{align*}
    \varphi(\alpha, \beta, \tau, \om)= &\alpha^2 \|P\|_F^2+ \beta^2 \|Q\|_F^2 +\tau^2 \|R\|_F^2+ (1-\om)^2 \|A\|_F^2 +(1-\om)^2 \|D\|_F^2 \\
    &+ 2 \alpha (\om -1) tr(PA)+ 2 \beta (\om -1) tr(QD)+ 2 (1-\om)^2 \|B\|_F^2\\
    &+ 2(1-\om)^2 \|C\|_F^2.
\end{align*}
%
Now we need to select the parameters $\alpha, \beta, \tau$ and $\om$ such that $\varphi(\alpha, \beta, \tau, \om)$ is very small. Since 
\begin{equation*}
     \lim_{\alpha,\, \beta,\, \tau,\, \om \rightarrow 0_{+}}\varphi(\alpha, \beta, \tau, \om) =  (1-\om)^2 \|A\|_F^2 +(1-\om)^2 \|D\|_F^2+ 2 (1-\om)^2 \|B\|_F^2+ 2(1-\om)^2 \|C\|_F^2,
\end{equation*}
 we can select $\om =1$ and $\alpha,\, \beta,\, \tau,\, \om \rightarrow 0_{+}$ such that $\varphi(\alpha, \beta, \tau, \om)\rightarrow 0_{+},$ and consequently, $\mathcal{N}_{\text{GSS}}\rightarrow \0.$
 
 In the sequel,  the GSS preconditioner is equivalently rewritten as: 
$ \P_{\text{GSS}}= \om \widetilde{\P}_{\text{GSS}},$
 where $$\widetilde{\P}_{\text{GSS}}=\bmatrix{\frac{\alpha}{\om} P+A& \0 & B^T\\ \0& \frac{\beta}{\om} Q+D & C\\ -B& -C^T &\frac{\tau}{\om}R}$$ can be regarded as a scaled preconditioner. A preconditioner is considered efficient if it closely approximates the coefficient matrix, and for given $\alpha, \beta, \tau >0,$ $\widetilde{\P}_{\text{GSS}}-\B=\frac{1}{\om}\Theta \rightarrow \0$ as $\om \rightarrow \infty.$ Thus, the GSS preconditioner shows enhanced efficiency for large values of $\om.$ Similar studies apply to the RGSS-I and RGSS-II preconditioners as well. The performance of the proposed preconditioners by varying the parameters is shown in the numerical experiment section.

  \section{Numerical experiments}\label{numerical}
  To demonstrate the effectiveness and robustness of the proposed preconditioners GSS, RGSS-I and RGSS-II over the existing ones within the Krylov subspace iterative method to solve the DSPP, in this section, we perform a few numerical experiments. We compare our proposed GSS, RGSS-I and RGSS-II PGMRES methods (abbreviated as ``GSS'', ``RGSS-I'' and ```RGSS-II'', respectively) with the GMRES (without preconditioning) method and PGMRES methods with  {block diagonal \cite{Susanne2023}, block preconditioner  \cite{ChangFeng2018}, dimension splitting  \cite{BenziDS2011}, relaxed dimension factorization \cite{BenziRDF2011}, and shift-splitting \cite{SSPDE2024}} preconditioners (abbreviated as ``BD'', ``BP''  ``DS'', ``RDF'', and ``SS'', respectively).  {Moreover, we have also compared proposed methods with BD preconditioner in conjunction with minimum residual method (MIRES) (abbreviated as ``BD-MINRES'').}  The numerical results are presented in terms of iteration counts (abbreviated as ``IT") and elapsed CPU time in seconds (abbreviated as ``CPU").

 The initial guess vector is set to ${\bm u}^{(0)}=\0\in \R^{n+l+m}$   for all iterative methods, and the method terminates if 
	$${\tt RES}:=\frac{\|\B{\bm u}^{(k+1)}-{\bm b}\|_2}{\|{\bm b}\|_2}<10^{-6}$$
 or if the maximum number of iterations exceeds $5000$. Numerical experiments are conducted in MATLAB R2024a on a Windows 11  system, using an Intel(R) Core(TM) i7-10700 CPU at 2.90 GHz with 16 GB of memory. 
 

 %
\begin{exam}\cite{Susanne2023, PDE-constrained2010}\label{exam1}
    \textbf{The Poisson control problem:} Consider the distributed control problem (PDE-constrained optimization problem) defined by:
    \begin{align}\label{pde1}
    \begin{array}{cc}
         & \underset{u,f}{\text{min}}~ \frac{1}{2} \|u-\hat{u}\|^2_{L_2(\Omega)}+\frac{ {\bm \nu}}{2} \|f\|^2_{L_2(\Omega)} \\
       \text{such that }  & -\nabla^2 u= f ~ \text{in} ~ \Omega,\\
       &\quad~~~u=g ~\text{on} ~ \partial\Omega,
    \end{array}
    \end{align}
  where $u$ is the state, and $\hat{u}$ is the desired state, $f$ is the control, $\Omega=[0,1]\times [0,1]$ is the domain with the boundary $\partial \Omega,$ and $0< { {\bm \nu}} \ll 1$  is the regularization parameter.  By discretizing \eqref{pde1} using the Galerkin finite element  method and then applying Lagrange multiplier techniques, we obtain the following system:
  \begin{align}\label{pde2}
      \bmatrix{{ {\bm \nu}}M &\0&K^T\\ \0&M&-M\\-K&M&\0}\bmatrix{u \\f \\ \bm{\lambda}}=\bmatrix{\0\\ b\\ -d },
  \end{align}
  where $M\in \R^{n\times n}$  and $K\in\R^{n\times n}$ are SPD mass matrix and discrete Laplacian, respectively. Note that by setting $ A =  {\bm \nu} M $, \( B = K \), \( C = -M \), and \( \mathbf{b} = [\mathbf{0}, b^T, -d^T]^T \),  \eqref{pde2} can be expressed in the form of the DSPP \eqref{SPP1}, where \( n = m = l \).

  The  MATLAB code downloaded from \cite{github}, generates the linear system of the form \eqref{pde2} by using the following setup: parameters in ``\textit{set\_def\_setup.m}'' are selected as \textit{def\_setup.bc}=`dirichlet', \textit{def\_setup.beta}=$1e-2$, \textit{def\_setup.ob}=1,  \textit{def\_setup.type}=`dist2d' and \textit{def\_setup.pow} = $5$, $6$ and $7$.  {For these selection of \textit{def\_setup.pow}, size of the coefficient matrix $\B$ of the DSPP \eqref{SPP1} is $2883,$ $11907$ and $48378,$ respectively.}
  	\begin{table}[ht!]
					\centering
			\caption{Experimental results of  GMRES,  BD,  {BD-MINRES, BP}, DS, RDF, SS, GSS,  {IGSS}, RGSS-I,  {IRGSS-I}, RGSS-II, and  {IRGSS-II}  PGMRES methods for Example \ref{exam1} when ${ {\bm \nu}}=0.1$}
			\label{tab1}
			\resizebox{12cm}{!}{
				\begin{tabular}{ccccc}
					\toprule
					Process& 	&def$\_$setup.pow $=5$ & def_setup.pow $=6$ &def_setup.pow $=7$  \\
					\midrule 
					&  size($\B$) &	  $2883 \times 2883$ & $11907\times 11907$&  $48378 \times 48378 $     \\
					\midrule
					\multirow{2}{*}{GMRES}& IT& $579$	& $2149$ & $\bf{--}$\\
					& CPU& $2.9618$&$367.9274$ & $\bf{--}$ \\	
					\midrule
                        \multirow{2}{*}{BD}& IT& $10$ &$10$ &$10$\\
					& CPU& $2.1418$&	$40.6138$ & $1512.0707$\\
                    \midrule
                         \multirow{2}{*}{BD-MINRES}&  IT& 	$9$	 & $8$ & $8$\\
					&  CPU&  $1.5042$&	 $28.7551$ &  $1122.2980$\\	
                    \midrule
                      \multirow{2}{*}{ BP}&  IT&  $3$ &  $3$ &  $3$\\
					&  CPU&  $0.7222
$ &  $17.1404$  &  $555.3527$\\	
					\midrule
					\multirow{2}{*}{DS}& IT& $31$& $40$& $51$	\\
					& CPU& $2.8514$ & $54.3105$ & $1706.9543$\\
					\midrule
					\multirow{2}{*}{RDF}& IT & $6$ &	$6$ & $8$\\
					& CPU& $0.8733$& $8.1585$ & $302.9091$ \\
                         \midrule
					\multirow{2}{*}{SS}& IT & $16$	& $22$ & $46$  \\
					& CPU  & $2.12460$ &	$45.9988$ &$3985.6141$\\
     \midrule
     \multirow{1}{*}{GSS}& IT& $2$ & $2$ &$2$ \\
				$\om_{exp} =30$	& CPU& $0.6909$ &$7.3233$ & $226.8625$\\
                         \midrule
                         \multirow{1}{*}{ IGSS}&  IT&  $2$ &  $2$ &  $2$ \\
				 $\om_{exp} =30$	&  CPU&  $0.6148$ & $6.3220$ &  $186.5088$\\
                         \midrule
                         \multirow{1}{*}{RGSS-I}& IT& $2$& $2$ &$2$\\
				$\om_{exp} =25$	& CPU& $0.7796$ & $7.72871$ & $249.7657$\\
                \midrule
                         \multirow{1}{*}{ IRGSS-I}&  IT&  $2$ &  $2$ &  $2$ \\
				 $\om_{exp} =25$	&  CPU&  $0.6536$ & $6.6898$ &  $179.4396$\\
                         \midrule
                         \multirow{1}{*}{RGSS-II}& IT& $2$ &$2$ &$2$\\
				$\om_{exp}=30$	& CPU& $0.7946$ & $6.8447$& $238.5576$\\
                \midrule
                         \multirow{1}{*}{ IRGSS-II}&  IT&  $2$ &  $2$ &  $2$ \\
				 $\om_{exp} =30$	&  CPU&  $0.5948$ & $6.1503$ &  $185.7962$\\
					\bottomrule
					\multicolumn{5}{l}{ $\bf{--}$ signifies that the method does not converge within $5000$ IT.}
				\end{tabular}
			}
		\end{table} 
  \begin{table}[ht!]
					\centering
			\caption{Experimental results of  GMRES,  BD,  {BD-MINRES, BP}, DS, RDF, SS, GSS,  {IGSS}, RGSS-I,  {IRGSS-I}, RGSS-II, and  {IRGSS-II}  PGMRES methods for Example \ref{exam1} when ${ {\bm \nu}}=0.001$}
			\label{tab2}
			\resizebox{12cm}{!}{
				\begin{tabular}{ccccc}
					\toprule
					Process&	&\textit{def\_setup.pow} $=5$ & \textit{def\_setup.pow} $=6$ & \textit{def\_setup.pow} $=7$  \\
					\midrule 
					&  size($\B$)	&  $2883 \times 2883$ & $11907\times 11907$&  $48378 \times 48378 $   \\
					\midrule
					\multirow{2}{*}{ GMRES}& IT& $919$ &$2254$	& $\bf{--}$\\
					& CPU& $5.8702$ & $409.5427$ & $\bf{--}$\\	
					\midrule
                      \multirow{2}{*}{BD}& IT& $19$ & $19$ & $19$\\
					& CPU& $4.3498$ & $84.9364$  & $3034.4227$\\	
                    \midrule
                      \multirow{2}{*}{ BP}&  IT&  $6$ &  $6$ &  $6$\\
					&  CPU&  $1.7077
$ &  $29.1209$  &  $965.1344$\\	
                    \midrule
                         \multirow{2}{*}{BD-MINRES}&  IT& 	$19$	 & $13$ & $13$\\
					&  CPU&  $2.8402$&	 $46.8062$ &  $1709.8236$\\
					\midrule
					\multirow{2}{*}{DS}& IT& $29$	& $39$ & $46$\\
					& CPU& $2.6072$ & $47.8576$ & $1603.9128$\\
					\midrule
					\multirow{2}{*}{RDF}& IT&$14$  &$11$ & $8$\\
					& CPU& $1.4320$ & $	14.3846$ & $324.8007$\\
                         \midrule
					\multirow{2}{*}{SS}& IT &$27$& $38$ & $68$	   \\
					& CPU &$3.3960$& $79.7055$ & $5481.9674$\\
     \midrule
     \multirow{1}{*}{GSS}& IT&  $2$ & $2$ & $2$\\
			$\om_{exp} =30$		& CPU&    $0.7159$ & $7.4041$ & $242.8123$  \\
                         \midrule
                         \multirow{1}{*}{ IGSS}&  IT&  $2$ &  $2$ &  $2$ \\
				 $\om_{exp} =30$	&  CPU&  $0.7884$ & $6.1971$ &  $182.3670$\\
                  \midrule
                         \multirow{1}{*}{RGSS-I}& IT& $2$ & $2$ & $2$\\
				$\om_{exp} =30$	& CPU& $0.78675$& $7.6515$ & $221.0918$ \\
                \midrule
                         \multirow{1}{*}{ IRGSS-I}&  IT&  $2$ &  $2$ &  $2$ \\
                          $\om_{exp} =25$	&  CPU&  $0.6358$ & $5.9984$ &  $188.0611$\\
                         \midrule
                         \multirow{1}{*}{RGSS-II}& IT&  $2$ & $2$ & $2$\\
				$\om_{exp} = 26$	& CPU& $0.6828$ & $6.9229$ & $247.4777$ \\
                \midrule
                 \multirow{1}{*}{ IRGSS-II}&  IT&  $2$ &  $2$ &  $2$ \\
				 $\om_{exp} =30$	&  CPU&  $0.6123$ & $6.1737$ &  $178.3082$\\
					\bottomrule
					\multicolumn{5}{l}{ $\bf{--}$ signifies that the method does not converge within $5000$ IT.}
				\end{tabular}
			}
		\end{table} 
  	\begin{figure}[ht!]
				\centering
				\begin{subfigure}[b]{0.32\textwidth}
					\centering
					\includegraphics[width=\textwidth]{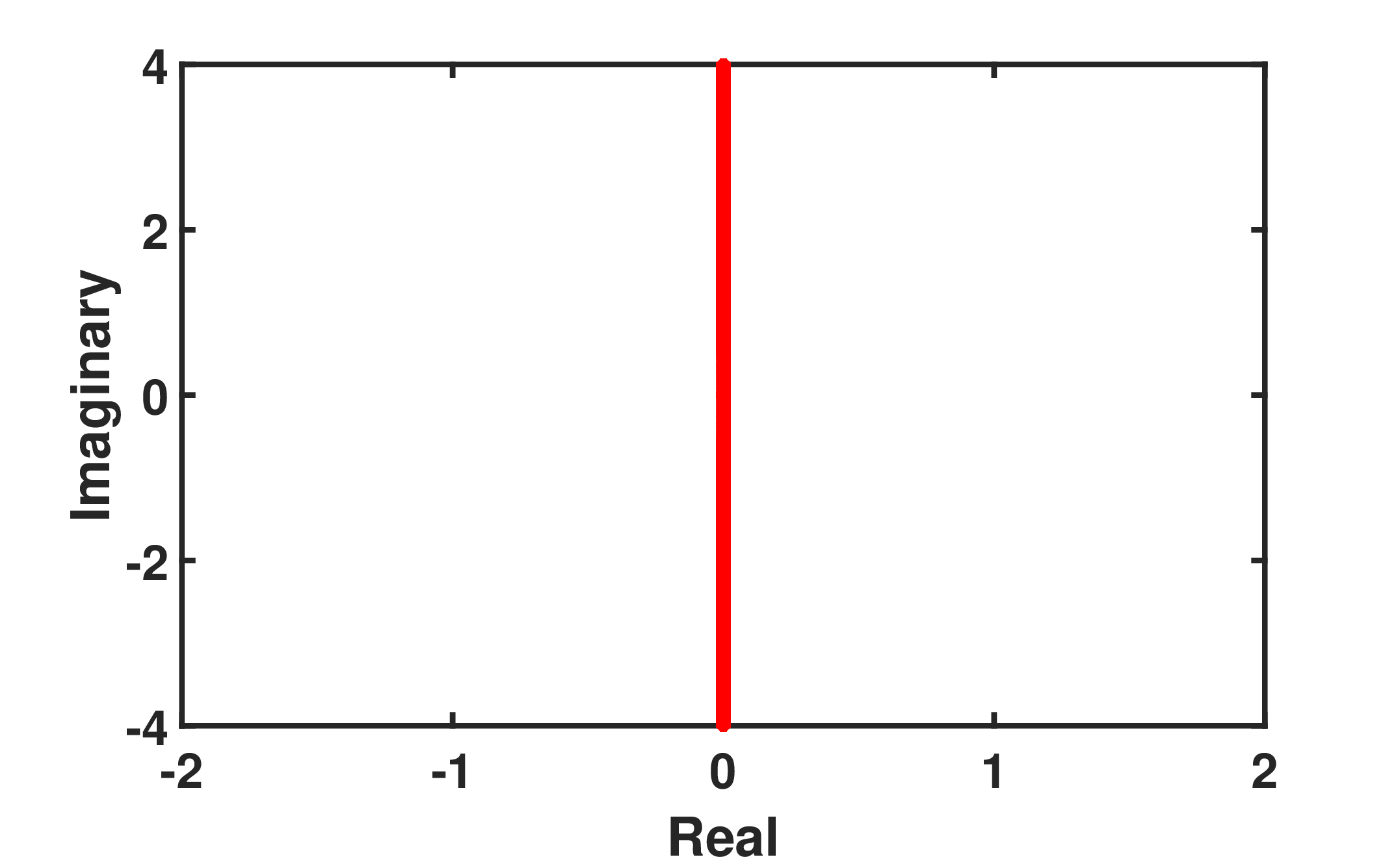}
					\caption{  $\B$ }
					\label{fig:original}
				\end{subfigure}
				\begin{subfigure}[b]{0.3\textwidth}
					\centering
					\includegraphics[width=\textwidth]{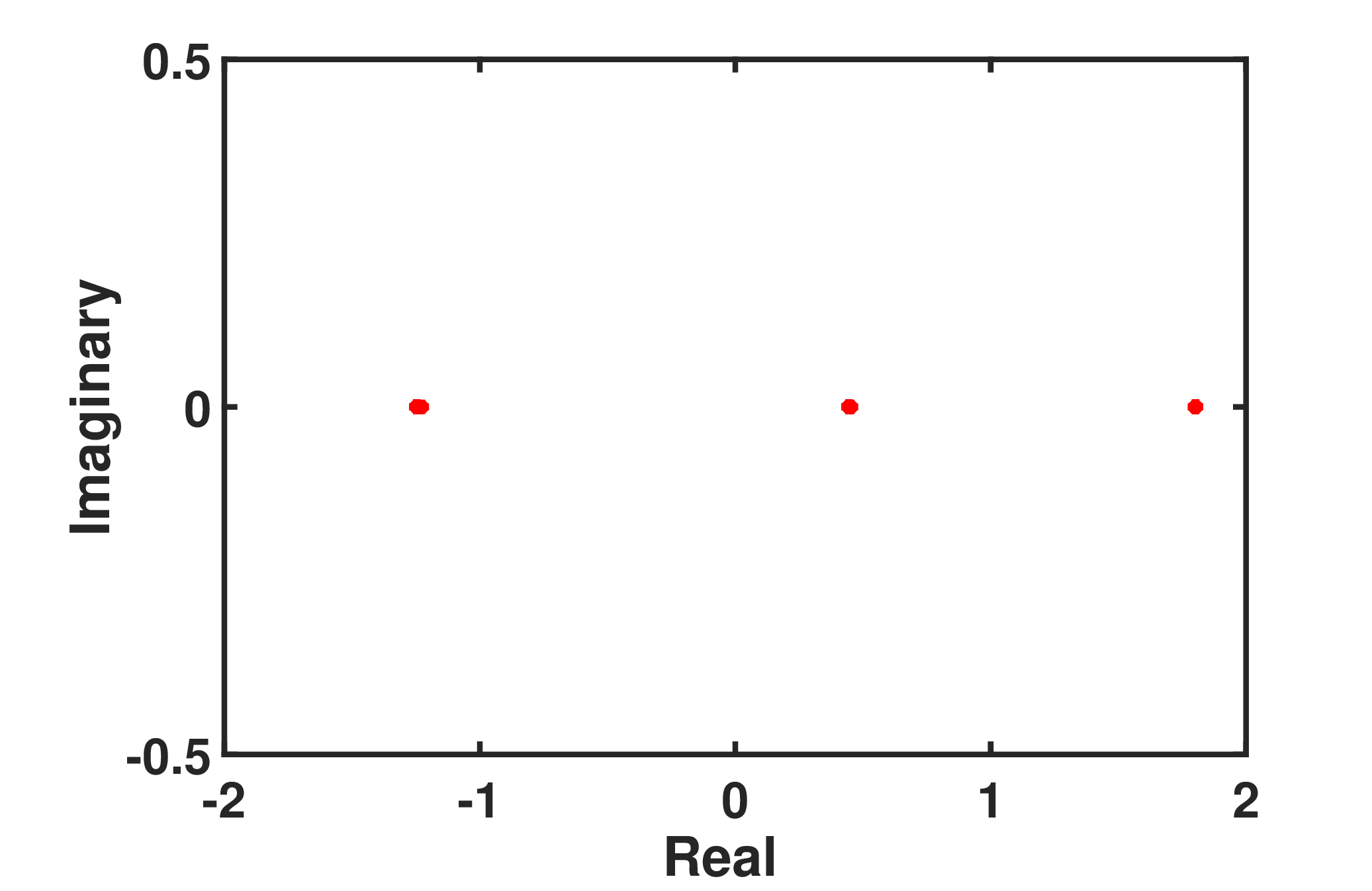}
					\caption{ $ \P_{\text{BD}}^{-1}\B$}
					\label{fig:BD}
				\end{subfigure}
                \begin{subfigure}[b]{0.3\textwidth}
					\centering
					\includegraphics[width=\textwidth]{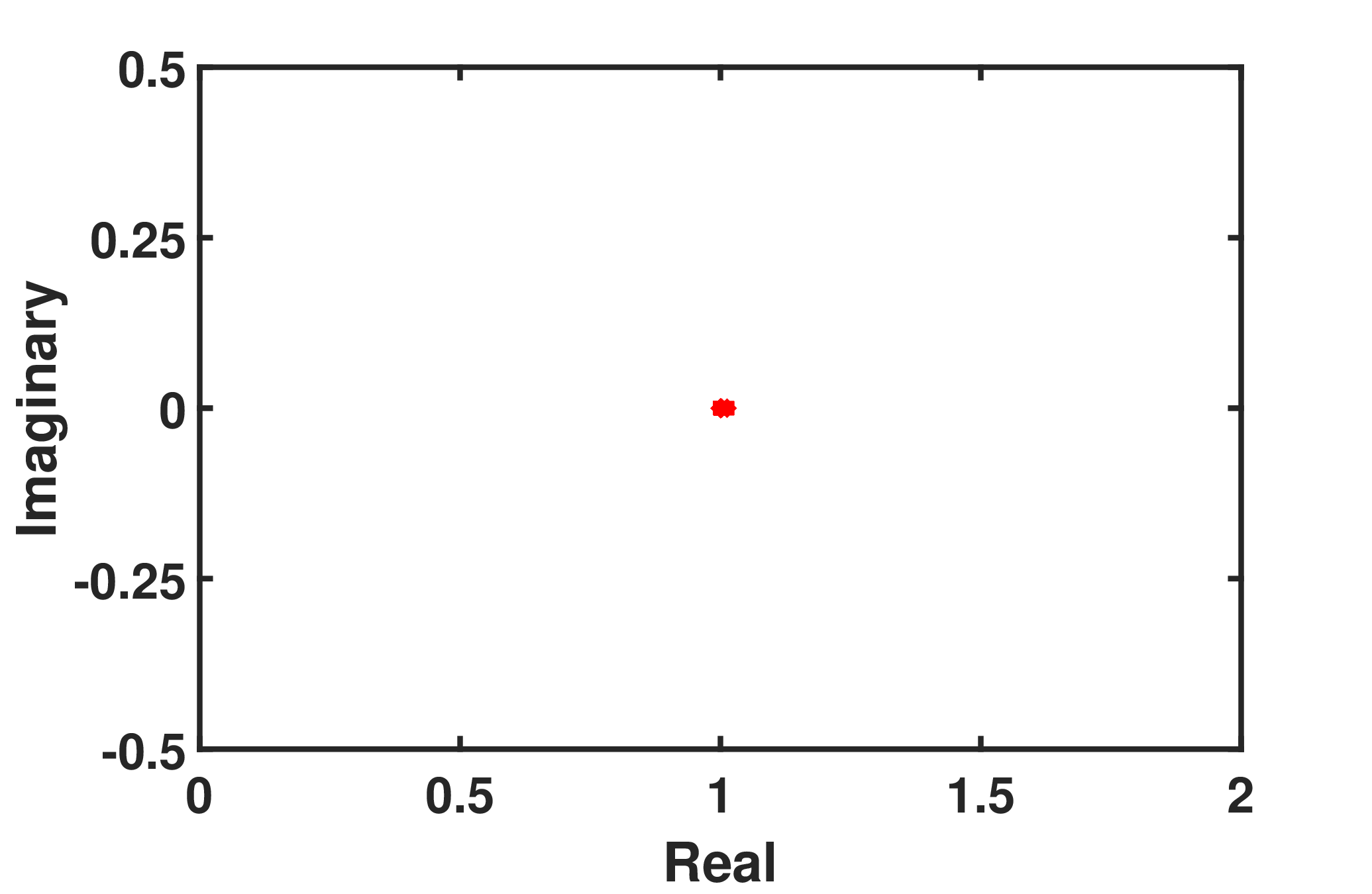}
					\caption{ $ \P_{\text{BP}}^{-1}\B$}
					\label{fig:BT}
				\end{subfigure}
    \begin{subfigure}[b]{0.3\textwidth}
					\centering
					\includegraphics[width=\textwidth]{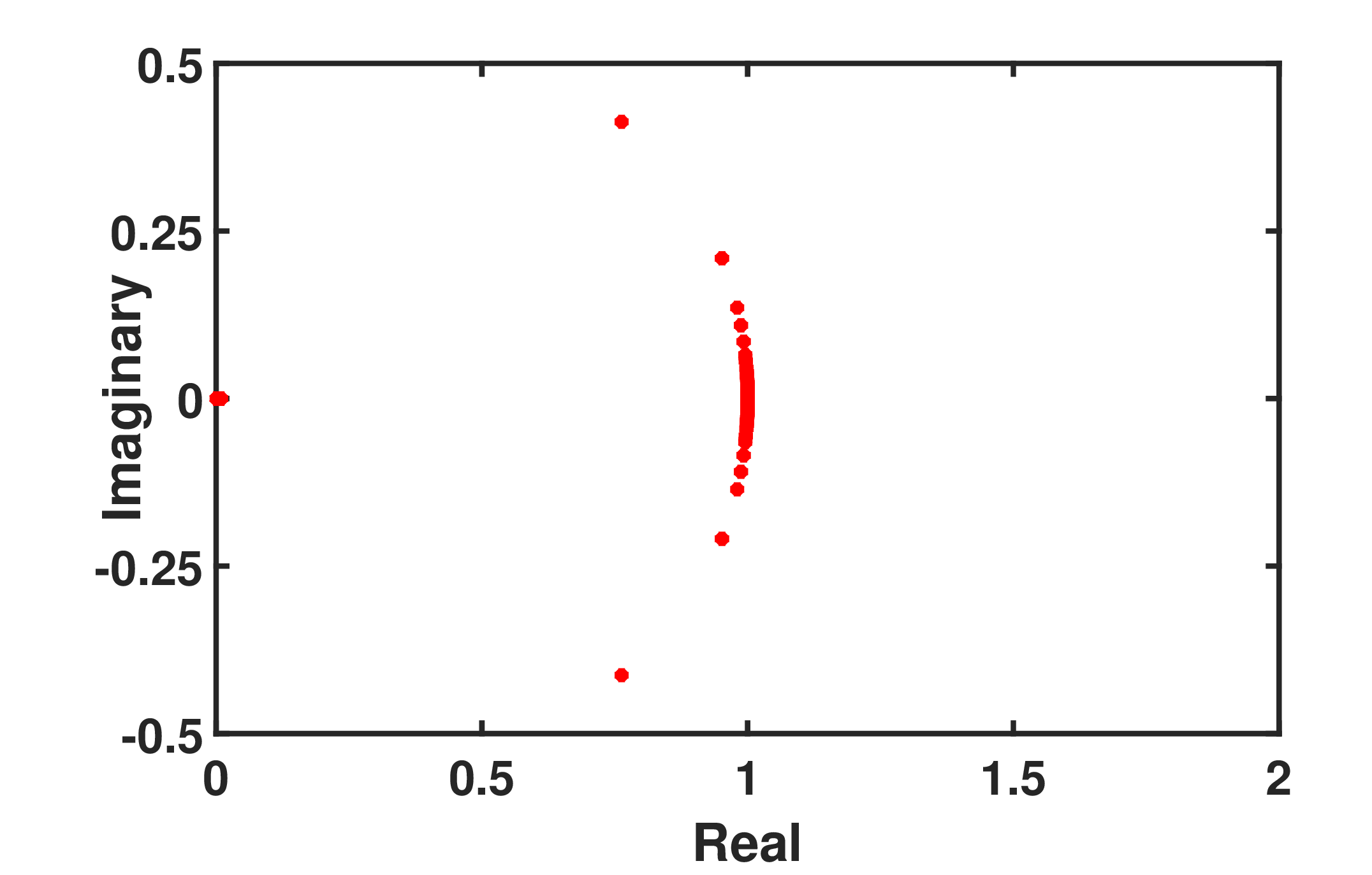}
					\caption{ $ \P_{\text{DS}}^{-1}\B$}
					\label{fig:DS}
				\end{subfigure}
    \begin{subfigure}[b]{0.33\textwidth}
					\centering
					\includegraphics[width=\textwidth]{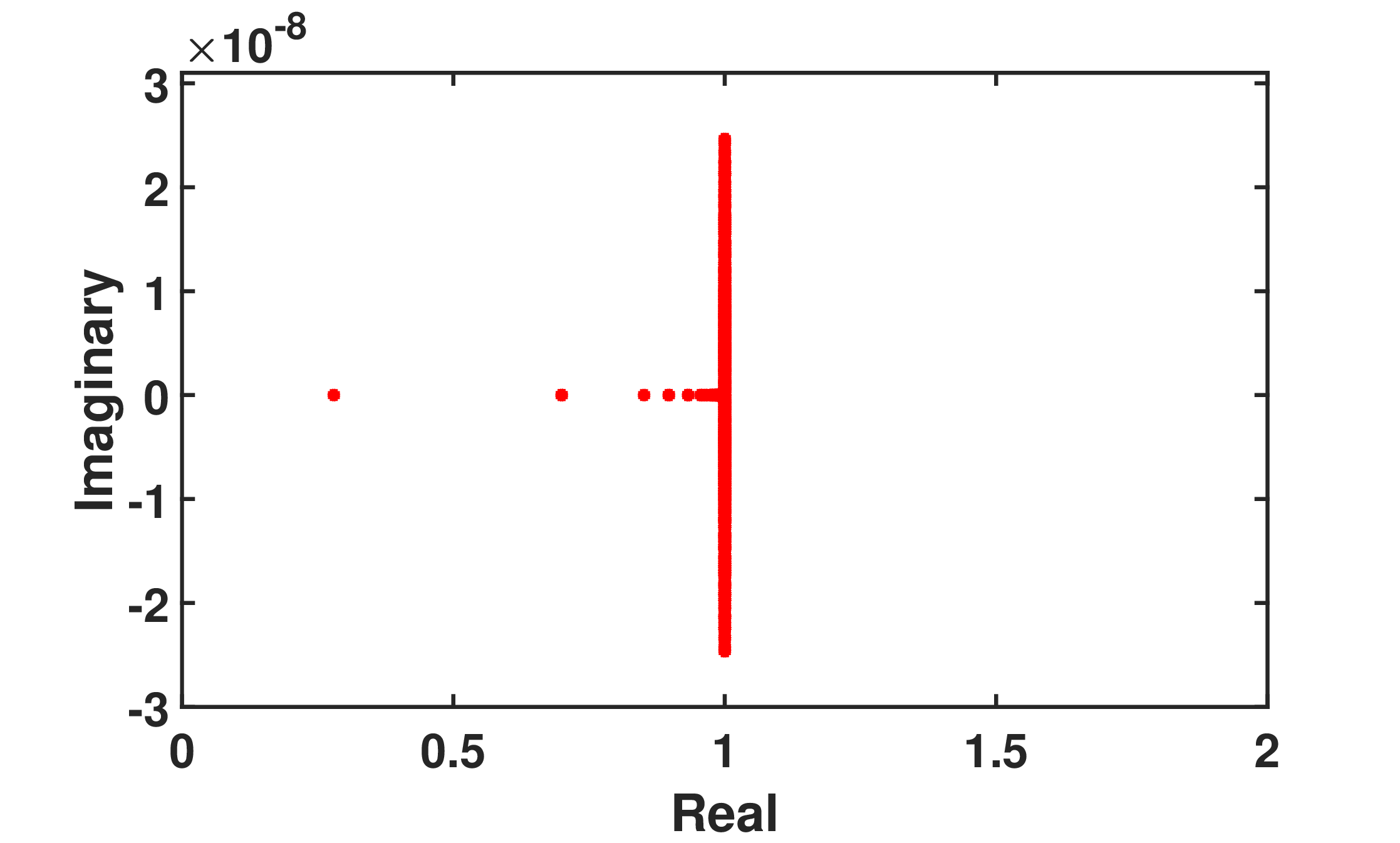}
					\caption{ $ \P_{\text{RDF}}^{-1}\B$}
					\label{fig:RDF}
				\end{subfigure}
    \begin{subfigure}[b]{0.31\textwidth}
					\centering
					\includegraphics[width=\textwidth]{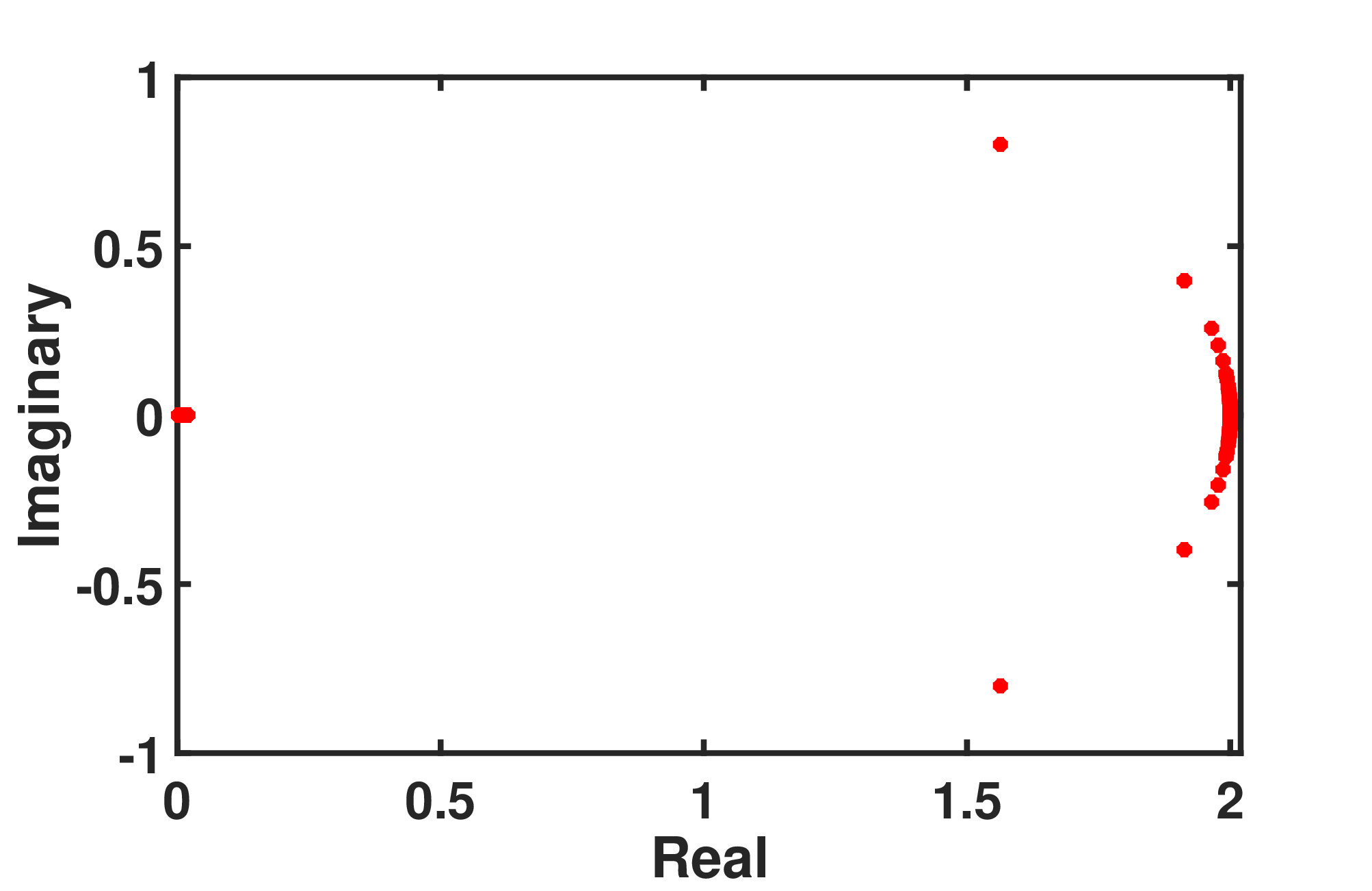}
					\caption{ $ \P_{\text{SS}}^{-1}\B$}
					\label{fig:SS}
				\end{subfigure}
				\begin{subfigure}[b]{0.32\textwidth}
					\centering
					\includegraphics[width=\textwidth]{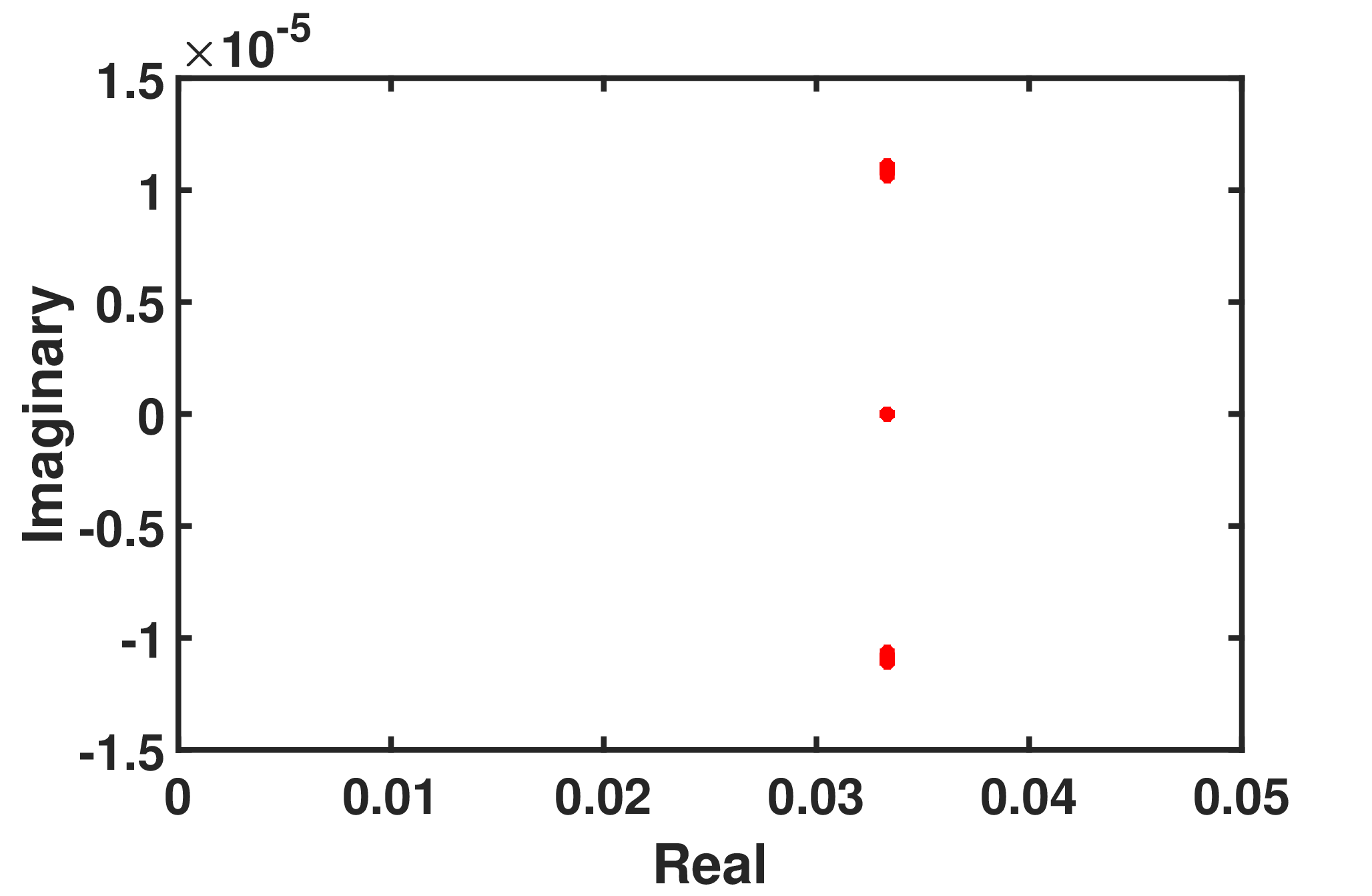}
					\caption{ $ \P_{\mathrm{GSS}}^{-1}\B$  $(\om=30)$ }
					\label{fig:InBD}
				\end{subfigure}
    \begin{subfigure}[b]{0.31\textwidth}
					\centering
					\includegraphics[width=\textwidth]{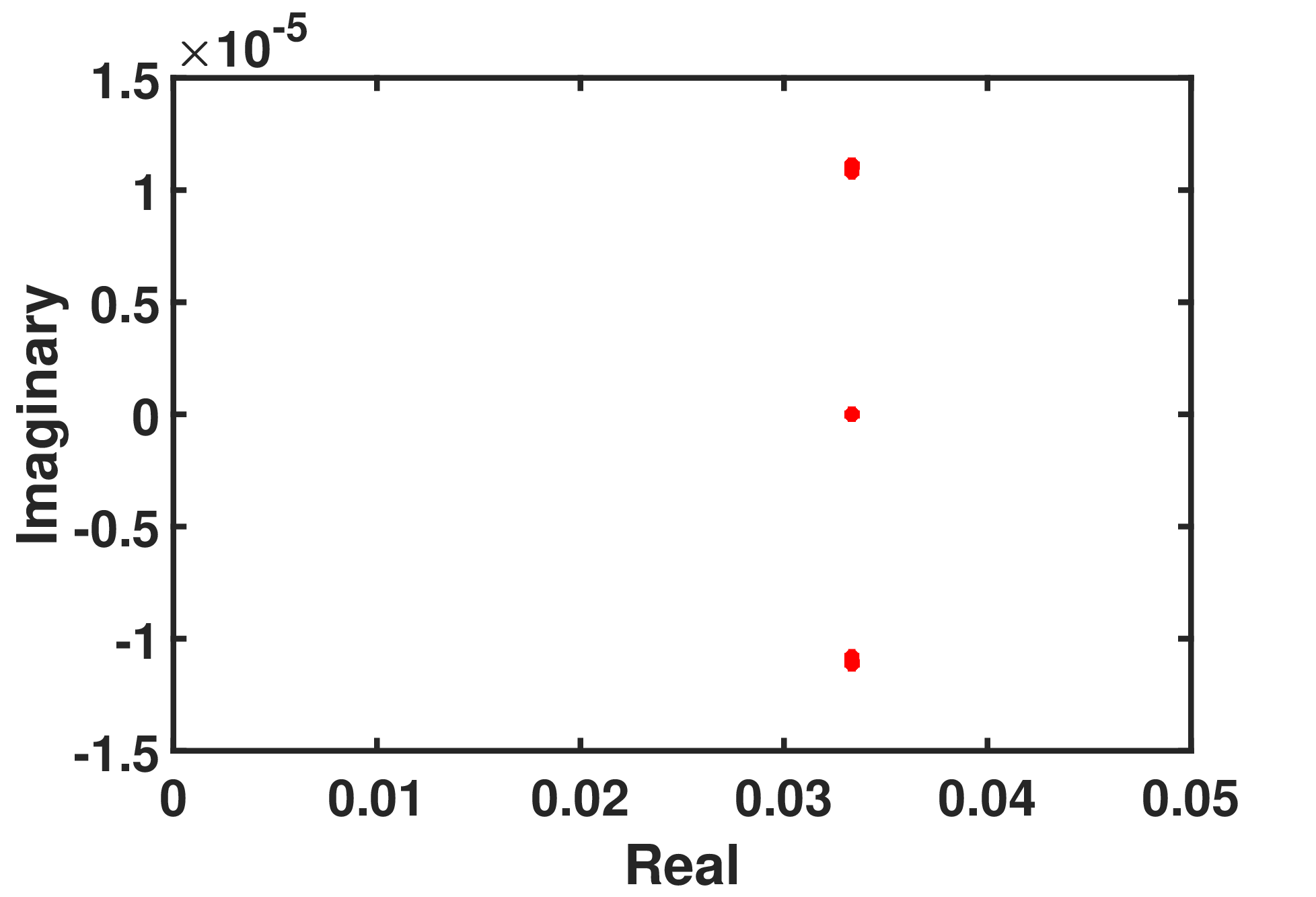}
					\caption{   $ \P_{\text{RGSS-I}}^{-1}\B$ $(\om=30)$ }
					\label{fig:RPSS_I}
				\end{subfigure}
     \begin{subfigure}[b]{0.31\textwidth}
					\centering
					\includegraphics[width=\textwidth]{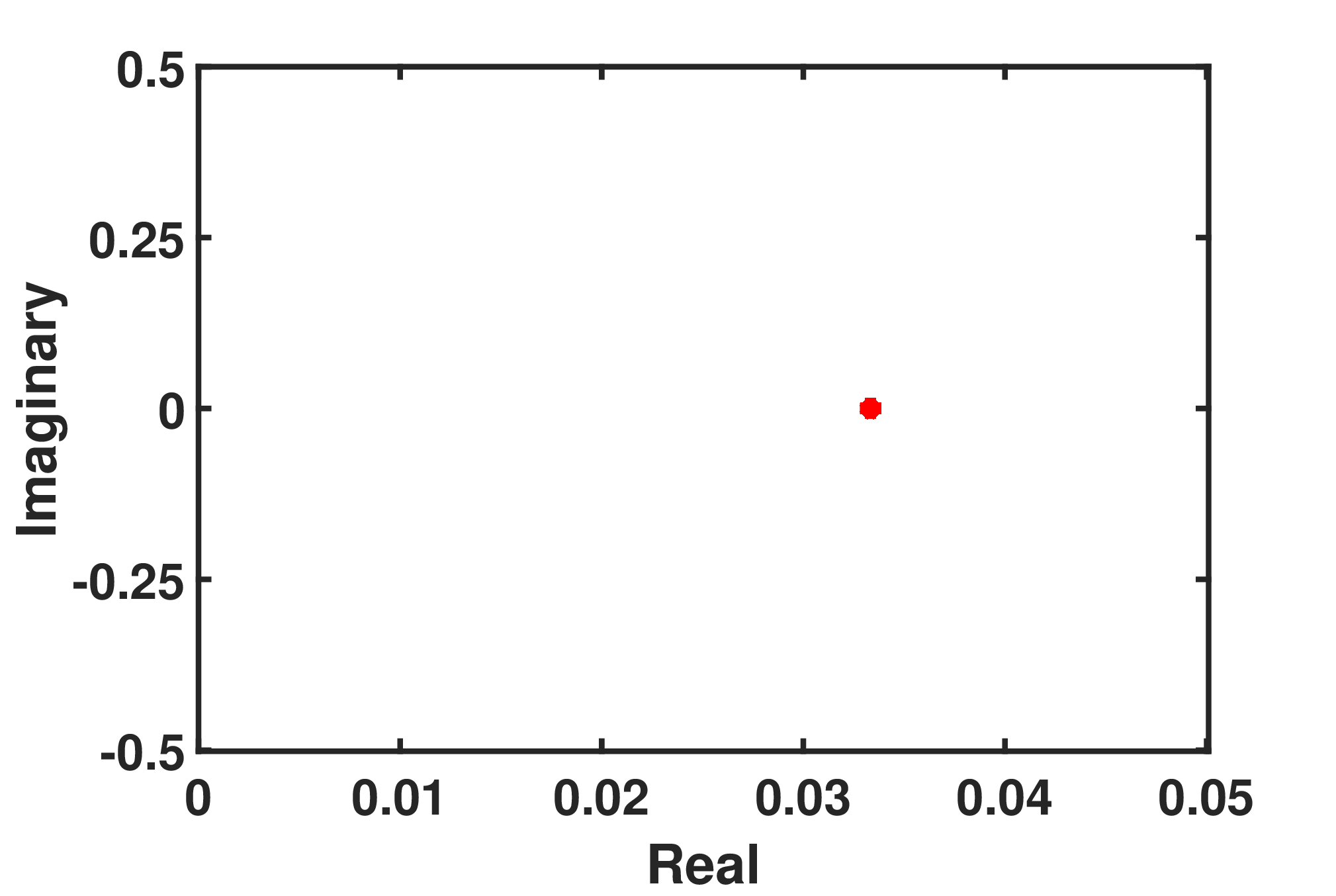}
					\caption{  $ \P_{\text{RGSS-II}}^{-1}\B$  $(\om=26)$}
					\label{fig:RGSS_II}
				\end{subfigure}
				\caption{ Eigenvalue distributions of $\B,$ $ \P_{\mathrm{BD}}^{-1}\B,$ $ \P_{\mathrm{BP}}^{-1}\B,$ $ \P_{\mathrm{DS}}^{-1}\B,$ $  \P_{\mathrm{SS}}^{-1}\B,$ $ \P_{\mathrm{GSS}}^{-1}\B,$ $  \P_{\text{RGSS-I}}^{-1}\B$ and $ \P_{\text{RGSS-II}}^{-1}\B$  for def$\_$set.pow $=5$ with $ {\bm \nu}=0.1$ for Example \ref{exam1}.}
				\label{fig1}
			\end{figure}

     {
    The BD preconditioner is utilized in its exact form, as specified in \cite[equation $(5.6)$]{Susanne2023}. We adopt the BP preconditioner as defined in \cite[equation $(2.1)$]{ChangFeng2018}. The DS preconditioner requires solving two linear subsystems with coefficient matrices \((\alpha I + A + \frac{1}{\alpha} B^T B)\) and \((\alpha I + A + \frac{1}{\alpha} C C^T)\). Similarly, the RDF preconditioner involves two linear subsystems with coefficient matrices \((A + \frac{1}{\alpha} B^T B)\) and \((A + \frac{1}{\alpha} C C^T)\). The SS preconditioner necessitates solving linear subsystems with coefficient matrices \((\alpha I + 2 \bm{\nu} M)\), \((\alpha I + M)\), and \((\alpha I + M(\alpha I + 2 \bm{\nu} M)^{-1})M + K(\alpha I + M)K^T\). For all these subsystems, we employ exact solvers, such as Cholesky factorization. Furthermore, in Algorithms 1, 2, and 3, all linear subsystems are solved using Cholesky factorization. For the inexact versions of our proposed preconditioners $\widetilde{\P}_{\GSS},$ $\widetilde{\P}_{\text{RGSS-I}}$ and $\widetilde{\P}_{\text{RGSS-II}}$ (abbreviated as ``IGSS", ``IRGSS-I" and ``IRGSS-II", respectively), we construct $\widetilde{P},\, \widetilde{P}_1,\, \widetilde{Q}, $ and $ \widetilde{Q}_1$ as follows:
consider $\widetilde{A}=L_AL^T_A, \,\widetilde{A}_p=L_{A_{p}}L^T_{A_{p}}, \widetilde{D}=L_DL^T_D,\, \widetilde{D}_q=L_{D_{q}}L^T_{D_{q}},$ $\widetilde{P}=\diag(B \widetilde{A}_{p}^{-1}B^T),$ $\widetilde{P}_1=\diag(B \widetilde{A}^{-1}B^T),$ $\widetilde{Q}=\diag(C^T\widetilde{D}_q^{-1}C),$ and $\widetilde{Q}_1=\diag(C^T\widetilde{D}^{-1}C),$ where $L_A,\, L_{A_p}\, L_D,$ and $L_{D_q}$ are the incomplete Cholesky factorization of $A,\, (\alpha P+\om A),\, D,$ and $(\beta Q+\om D),$ respectively,  produced by the Matlab function: 
$$ \textit{\textbf{ichol(A, struct(`type', `ict', `droptol', \text{1e-4}, `michol', `off'))}}.$$}

\vspace{-2mm}\noindent\textbf{Parameter selection:} For the DS preconditioner, we choose the parameter $\alpha$ (denoted by $\alpha_{\text{DS}}$) as  {follows} \cite{BenziDS2011}: $$\alpha_{\text{DS}}=\frac{\sqrt{tr(A^TA)+2tr(BB^T)}+\sqrt{tr(D^TD)+2tr(C^TC)}}{2(n+m+l)}.$$ For the RDF preconditioner, the parameter \(\alpha\) is selected from the interval \((0,1)\) with a step size of 0.01. The optimal performance, in terms of minimal CPU times, is achieved with smaller values of \(\alpha\) as found in  \cite{BenziRDF2011}. For the SS preconditioner, we take $\alpha=0.01.$
   For the GSS, RGSS-I, and RGSS-II preconditioners, as well as for their inexact forms, the parameters are chosen as follows: $\alpha= \beta=0.01,$ $\tau=0.001,$ $P=A,$ $Q=CC^T,$ and $R=I.$ The optimal parameter $\om$ (denoted by $\om_{exp}$) is determined experimentally within the interval \([2,30]\) with step size one, which yields minimal CPU times. 
   
\vspace{2mm}
\noindent \textbf{Numerical results:}   The numerical results for GMRES and various PGMRES methods with $ {\bm \nu}=0.1$ and $0.001$ are presented in Tables \ref{tab1} and \ref{tab2}. We observe that the GMRES method exhibits a significantly slower convergence rate compared to all other PGMRES methods, even does not converge within $5000$ iterations when \textit{def\_setup.pow} $= 7$. On the other hand, we observe that our proposed preconditioners outperform all the compared preconditioners in terms of both IT and CPU times. For the DS and SS preconditioners, the IT increases as the size of the saddle point matrix $\B$ grows. Whereas the proposed GSS, RGSS-I and RGSS-II preconditioners maintain a consistent IT regardless of matrix size.  {Moreover, from Tables \ref{tab1} and \ref{tab2}, we observe that the inexact preconditioners IGSS, IRGSS-I, and IRGSS-II require less CPU time to converge, demonstrating more efficiency compared to their exact counterparts GSS, RGSS-I, and RGSS-II, respectively.}

\vspace{2mm}
\noindent \textbf{Eigenvalue distributions:} In order to better illustrate the effectiveness of the proposed GSS, RGSS-I and RGSS-II preconditioners, the eigenvalue distribution of $\B$ and preconditioned matrices $ \P_{\mathrm{BD}}^{-1}\B,$ $ \P_{\mathrm{BP}}^{-1}\B,$ $ \P_{\mathrm{DS}}^{-1}\B,$ $ \P_{\mathrm{RDF}}^{-1}\B,$ $ \P_{\mathrm{SS}}^{-1}\B,$ $ \P_{\text{GSS}}^{-1}\B\, {(\om=30)},$ $ \P_{\text{RGSS-I}}^{-1}\B\, {(\om=30)}$ and $ \P_{\text{RGSS-II}}^{-1}\B\, {(\om=26)}$ (with ${\bm \nu}=0.1$) are displayed in Figure \ref{fig1}. From Figure \ref{fig1}, we observe that the eigenvalues of the preconditioned matrices $ \P_{\text{GSS}}^{-1}\B,$ $ \P_{\text{RGSS-I}}^{-1}\B,$ and $ \P_{\text{RGSS-II}}^{-1}\B$ have clustered better than the coefficient matrix $\B,$ and preconditioned matrices $ \P_{\mathrm{BD}}^{-1}\B,$ $ \P_{\mathrm{BP}}^{-1}\B,$ $ \P_{\mathrm{DS}}^{-1}\B,$ $ \P_{\mathrm{RDF}}^{-1}\B$ and $ \P_{\mathrm{SS}}^{-1}\B.$  {The spectral radius of the saddle point matrix $\B$ is $\vartheta(\B)=3.9872$. The eigenvalues of \(\P_{\mathrm{BD}}^{-1} \B\) are all real with $\vartheta(\P_{\mathrm{BD}}^{-1}\B)= 1.8019$ and  $\vartheta({\P}_{\text{BP}}^{-1}\B)=1.0128.$
 The preconditioned matrices \(\P_{\mathrm{DS}}^{-1} \B\) and \(\P_{\mathrm{RDF}}^{-1} \B\) both have a spectral radius of 1. However, \(\P_{\mathrm{DS}}^{-1} \B\) has non-real eigenvalues with imaginary parts of order \(10^{-1}\) or \(10^{-2}\), whereas the imaginary parts of the non-real eigenvalues of \(\P_{\mathrm{RDF}}^{-1} \B\) are significantly smaller, on the order of \(10^{-8}\). The computed spectral radius for   preconditioned matrices $ \P_{\mathrm{SS}}^{-1}\B,$ $ \P_{\text{GSS}}^{-1}\B,$ $ \P_{\text{RGSS-I}}^{-1}\B$ and $ \P_{\text{RGSS-II}}^{-1}\B$ is $\vartheta(\P_{\mathrm{SS}}^{-1}\B)= 1.9995$, $\vartheta(\P_{\GSS}^{-1}\B)= 0.0333$, $\vartheta(\P_{\text{RGSS-I}}^{-1}\B)= 0.0333$, and $\vartheta(\P_{\text{RGSS-I}}^{-1}\B)=0.0385,$ respectively. However, the complex part of the non-real eigenvalues of the SS preconditioned matrix is in the range of $10^{-1},$ whereas for GSS and RGSS-I, they are in the range of $10^{-5}.$  There are at least $961 (=n)$ eigenvalues of RGSS-I preconditioned matrix is equal to $0.0333$ and for the RGSS-II preconditioned matrix, at least $1922 (=n+l)$ real eigenvalues are equal to $0.0385.$ Thus, we can conclude that GSS, RGSS-I, and RGSS-II preconditioned matrices have better clustered spectrums compared to other preconditioned matrices. This indicates enhanced computational efficiency, highlighting the effectiveness of the GSS, RGSS-I, and RGSS-II preconditioners.}

\vspace{2mm}
\noindent {\textbf{Influence of the parameters ${\alpha, \beta, \tau, \om}$:}}
 \begin{figure}[ht!]
				\centering
				\begin{subfigure}[b]{0.4\textwidth}
					\centering
					\includegraphics[width=\textwidth]{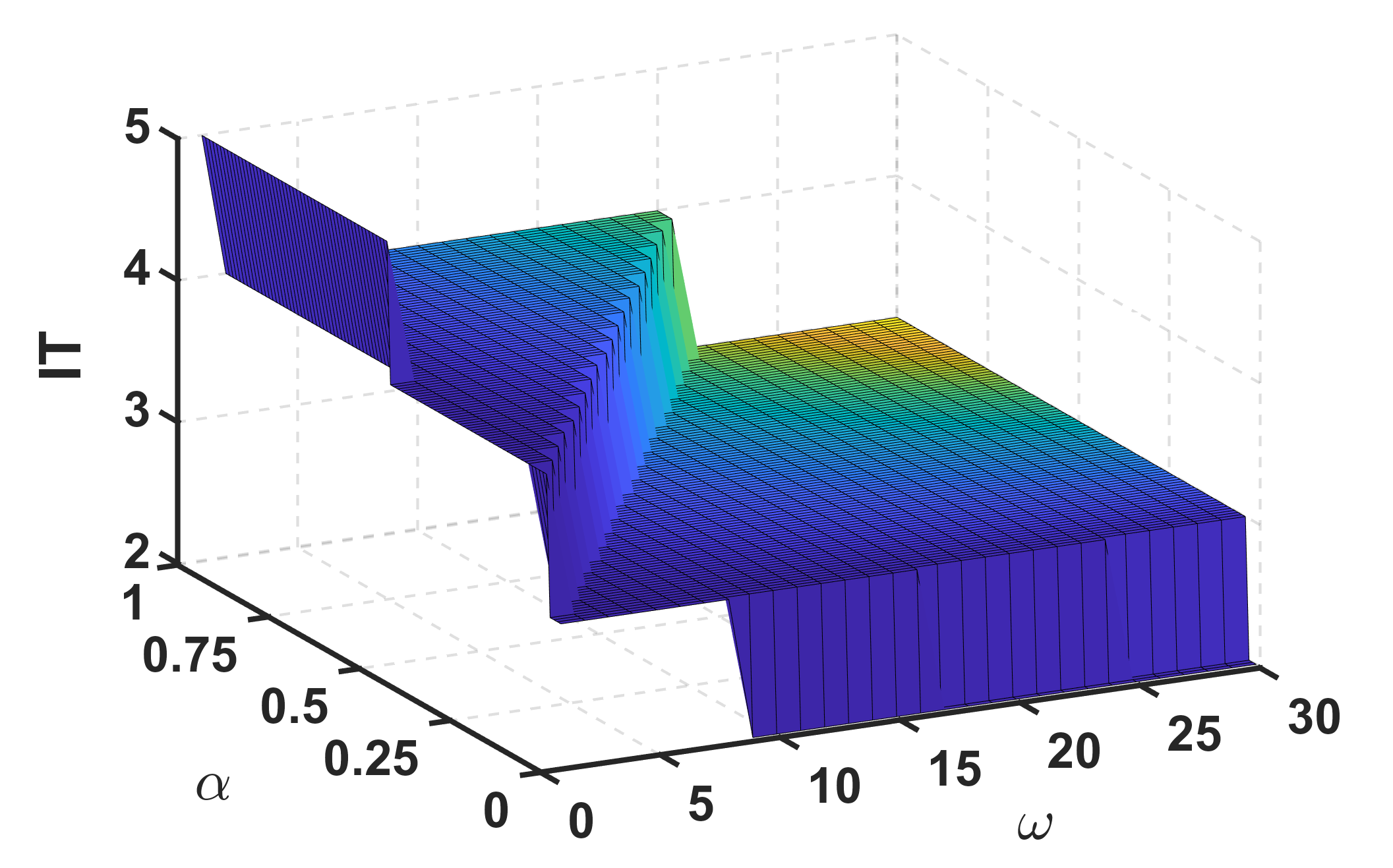}
					\caption{\footnotesize By varying \(\alpha = \beta\)  within the interval \([0.01, 1]\) and \(\omega\)  within the range \([1, 30]\) for the GSS preconditioner}
					\label{fig:GSS_AL}
				\end{subfigure} 
    \begin{subfigure}[b]{0.4\textwidth}
					\centering
					\includegraphics[width=\textwidth]{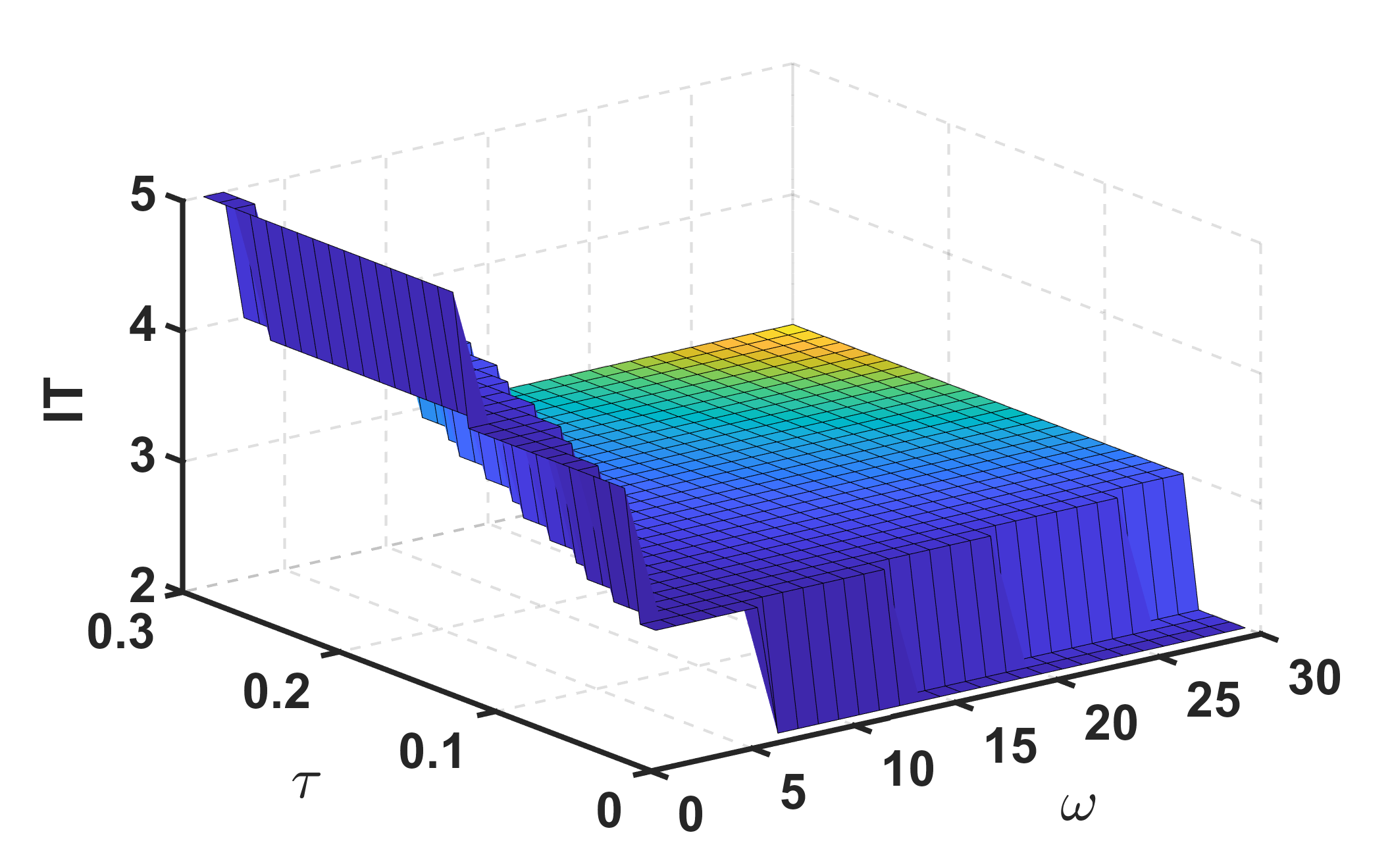}
					\caption{\footnotesize By varying $\tau$ within the interval $[0.01, 0.3]$ and $\om$ within the interval
     $[1,30]$ for the RGSS-II preconditioner}
					\label{fig:RGSS_II_tau}
				\end{subfigure} 
\caption{Convergence curves of the GSS and RGSS-II RGMRES methods varying the parameters $\alpha,$ $\beta,$ $\tau,$ $\om$   for Example \ref{exam1} with $ {\bm \nu}=0.1$.}
 \label{fig: parameters}
			\end{figure}
  To demonstrate the influence of the parameter on the performance of the proposed preconditioners, we present graphs of IT counts versus parameters for the GSS and RGSS-II preconditioners in Figure \ref{fig: parameters}. For the GSS preconditioner, we vary \(\alpha = \beta\) from $0.01$ to $1$ with a step size of $0.01$ and \(\om\) from $1$ to $30$ with a step size of one. For the RGSS-II preconditioner, we vary \(\tau\) from $0.01$ to $0.3$ with a step size of $0.01$ and \(\om\) from $1$ to $30$ with a step size of one. We can draw the following observation from Figure \ref{fig: parameters}:
  \begin{itemize}
      \item For both preconditioners, IT exhibits minimal sensitivity to variations in the parameters.
      \item Although, for small values of $\om$, IT increases when $\alpha$ increase for GSS preconditioner and $\tau$ increases for RGSS-II  preconditioner. Nonetheless, as $\om$ increases, IT decreases, even as the magnitude of $\alpha$ and $\tau$ continue to grow.
  \end{itemize}
  Therefore, the proposed preconditioners achieve high efficiency when \(\alpha\), \(\beta\) and \(\tau\) are kept small and \(\om\) is large.

\vspace{2mm}
\noindent	{\textbf{Condition number analysis:}} 
   \begin{figure}[ht!]
				\centering
				\begin{subfigure}[b]{0.5\textwidth}
					\centering
					\includegraphics[width=\textwidth]{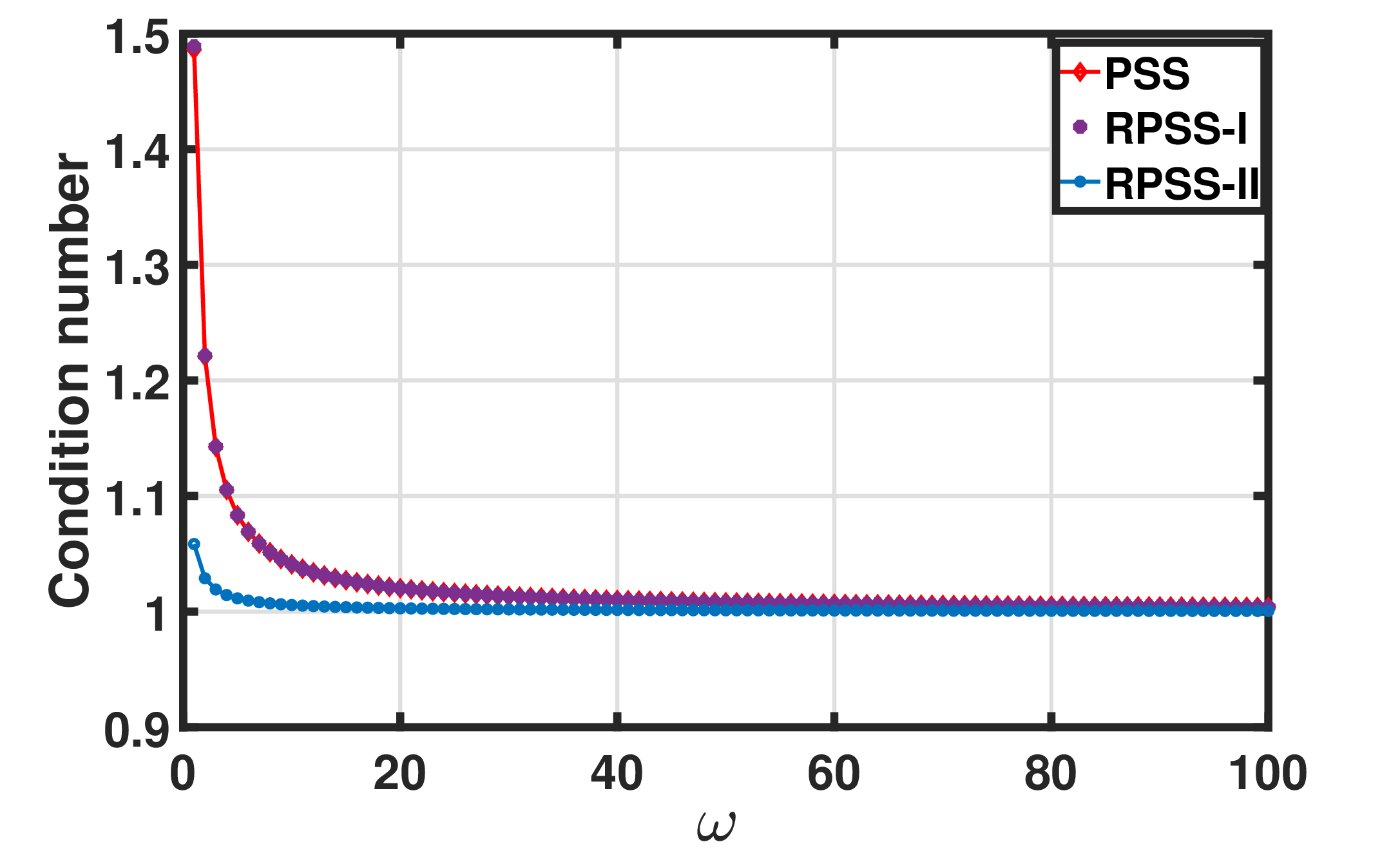}
				\end{subfigure} 
\caption{Relationship between condition numbers of the preconditioned matrices $ \P_{\text{GSS
}}^{-1}\B,  \P_{\text{RGSS-I}}^{-1}\B$ and $ \P_{\text{RGSS-II}}^{-1}\B$ varying the parameter $\om$ in $[1, 100]$  with $ {\bm \nu}=0.1$ for Example \ref{exam1}.}
 \label{fig:condition}
			\end{figure}
To evaluate the robustness of the proposed GSS, RGSS-I, and RGSS-II preconditioners, we assess the condition numbers of the preconditioned matrices $ \P_{\text{GSS
}}^{-1}\B,  \P_{\text{RGSS-I}}^{-1}\B$ and $ \P_{\text{RGSS-II}}^{-1}\B$. For any nonsingular matrix \(A\), the condition number \(\kappa(A)\) is defined as $
\kappa(A) := \|A^{-1}\|_2 \|A\|_2,$ 
where $\|A\|_2$ denotes the spectral norm of $A.$  {The condition number of $\B$ is $3.6396e+05$, which is moderately hight. On the other hand, in Figure \ref{fig:condition},} we depict the effect of the parameter $\om$ ranges from 1 to 100 with a step size of one on the preconditioned matrices $ \P_{\text{GSS
}}^{-1}\B,  \P_{\text{RGSS-I}}^{-1}\B$ and $ \P_{\text{RGSS-II}}^{-1}\B$  for the case \textit{def\_setup.pow} $=5$. We observe that for all values of \(\om\), the condition numbers of the preconditioned matrices remain within the range \([1, 1.5]\). This indicates that the preconditioned systems are well-conditioned, demonstrating that the GSS, RGSS-I, and RGSS-II preconditioners are robust and effective.

 \vspace{2mm}
\noindent\textbf{Analyzing the sensitivity of the solution:}
 \begin{figure}[ht!]
				\centering
				\begin{subfigure}[b]{0.5\textwidth}
					\centering
					\includegraphics[width=\textwidth]{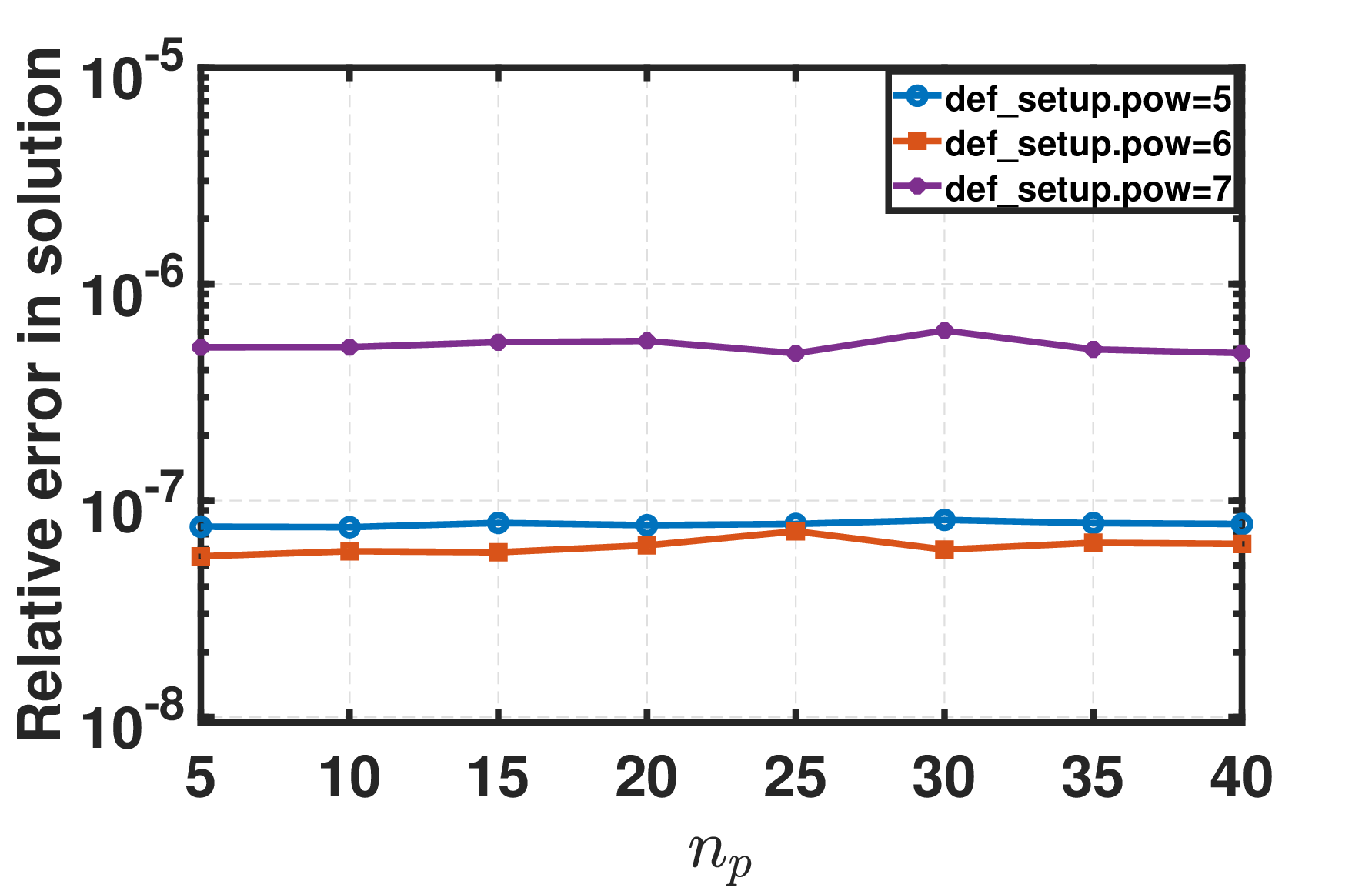}
				\end{subfigure} 
\caption{Behaviour of the relative errors in the solution with with increasing $n_p$ for \textit{def\_setup.pow} $=5,6,7$  with $ {\bm \nu}=0.1$ for Example \ref{exam1}.}
 \label{fig: sensitivity}
			\end{figure}
   To assess the reliability of the solution computed using the GSS PGMRES method when the coefficient matrix is slightly perturbed, we consider the perturbed system $(\B + \Delta \B) \bm{\widetilde{u}} = \bm{b}$ of the DSPP \eqref{SPP1}. We compute $\Delta \B$ by adding some Gaussian noise to  $B$ and $C$  in the following manner: 
$$\Delta B=\epsilon \cdot n_p\cdot std(B).*randn(m,n)\quad\text{and} \quad \Delta C=\epsilon \cdot n_p \cdot std(C).*randn(p,m),$$
     where  $n_p$ is the noise percentage,  $\epsilon=10^{-6},$ $randn(m,n)$ and $std(X)$ denote the random matrices of order $m\times n$ and the standard deviation of the matrix $X,$ respectively.  We conducted numerical tests for the GSS PGMRES method for \textit{def\_setup.pow} $=5,6,7$ with $\beta=0.1$. Let $\bm{u}$ and $\bm{\widetilde{u}}$ are the computed solutions of the DSPP and its perturbed system, respectively. We compute the relative error in the solution, given by $\dm{\frac{\|\bm{\widetilde{u}}-\bm{u}\|_2}{\|\bm{u}\|_2}}$ as the noise percentage $n_p$ increases from $5\%$ to $40\%$ with step size $5\%.$ These results are presented in Figure \ref{fig: sensitivity}. We observe that as $n_p$ increases, the relative error of the solution consistently remains below $10^{-6}$. This indicates that the solution \( \bm{u} \) is insensitive in nature to small perturbations in the coefficient matrix and demonstrates the robustness of the proposed GSS PGMRES method. Similar observations are also noticed for RGSS-I and RGSS-II PGMRES methods.
 
 \end{exam}
 { \begin{exam}\label{EX2}
    \textbf{The leaky lid driven cavity problem \cite{Matlab2007, CFD2005}:} We consider the DSPP arising from the discretization of the incompressible Stokes flow problem, which is given by
    \begin{align}\label{stokes}
				\begin{array}{rl}
				      -\Delta {\bm u}+\nabla \bm{v}=\0, \quad &~\text{in} \quad \Xi,\\ 
				\quad \quad \quad\nabla \cdot {\bm u}=0,  \quad &~\text{in}\quad \Xi,\\ 
                \quad \quad \quad \bm{u}=\0, \quad &~\text{on}\quad \Xi\backslash  \partial \Xi_{lid},\\ 
                \quad \quad \quad \bm{u}_x=1, \quad& ~\text{on}\quad   \partial \Xi_{lid},
                \end{array}
   			\end{align}
            where \(\Xi = (-1,1) \times (-1,1)\) represents the domain and \(\partial \Xi_{\text{lid}} = [-1,1] \times \{1\}\) denoting the lid boundary. The velocity field is given by the vector-valued function \(\bm{u} = (\bm{u}_x, \bm{u}_y)\), while the scalar function \(\bm{v}\) represents the pressure. Moreover, \(\Delta\) denotes the Laplacian operator in \(\mathbb{R}^2\), \(\nabla\) represents the gradient, and \(\nabla \cdot\) denotes the divergence.

            \begin{table}[ht!]
			 	\centering
				\caption{ Experimental results of  GMRES,  BD,  {BD-MINRES}, DS, RDF, SS, GSS, RGSS-I and RGSS-II  PGMRES methods for Example  \ref{EX2}.}
				\label{tab3}
				\resizebox{12cm}{!}{
					\begin{tabular}{ccccccc}
						\toprule
						Process& Grids	& $16\times 16$& $32\times 32$ & $64\times 64$ &  $128\times 128$\\
						\midrule 
						& size($\B$)	 &$770 \times 770$& $2946 \times 2946$ & $11522 \times 11522$& $45570 \times 45570$ \\
						\midrule
					\multirow{2}{*}{GMRES}& IT& $246$&$422$ &$692$	& $1169$\\
					& CPU& $0.2042$&$2.7629$ & $158.8824$ & $1880.7285$\\	
					\midrule
                      \multirow{2}{*}{BD}& IT& $38$& 	$38$	 & $38$ &$38$\\
					& CPU&$0.5025$& $5.2317$&	 $104.2237$  & $1067.3401$\\	
                    \midrule
                         \multirow{2}{*}{BD-MINRES}&  IT& $32$&$32$ & $32$ & $32$\\
					&  CPU&  $0.3213$& $4.5613$ & $90.6470$ & $650.6753$\\
					\midrule
					\multirow{2}{*}{DS}& IT&$30$& $39$	& $49$ & $59$\\
					& CPU&$0.5662$ &$3.6632$ & $104.0399$ & $885.1938$\\
					\midrule
					\multirow{2}{*}{RDF}& IT&$11$ &$13$  &$12$ & $11$\\
					& CPU& $0.4619$&$1.5057$ & $	26.8306$ & $501.8811$\\
                         \midrule
					\multirow{2}{*}{SS}& IT &$11$&$18$& $32$ & $53$	   \\
					& CPU &$0.3426$&$1.9709$& $48.6013$ & $1675.7429$\\
     \midrule
     \multirow{1}{*}{GSS}& IT& $2$&  $2$ & $2$ & $3$\\
			$\om_{exp} =25$		& CPU& $0.2040$& $0.7120$& $6.0955$ & $205.7726$  \\
                         \midrule
                         \multirow{1}{*}{RGSS-I}& IT&$2$&  $2$ & $2$ & $3$\\
				$\om_{exp} =29$	& CPU&$0.2296$& $0.6767$& $5.3381$ & $207.4075$ \\
                         \midrule
                         \multirow{1}{*}{RGSS-II}& IT& $2$ &$2$ & $2$ & $3$\\
				$\om_{exp} = 29$	& CPU&$0.2199$ & $0.6125$ & $5.4854$ & $209.7603$ \\
					\bottomrule
					\end{tabular}
				}
			\end{table} 

 Discretizing the systems in \eqref{stokes} using the IFISS software package \cite{Matlab2007} with $Q2-Q1$ finite elements on uniform leads to the DSPP \eqref{SPP1}. The right hand side vector $\bm{b}\in \R^{n+l+m}$ is chosen so that the exact solution of the DSPP is $\bm{u}=[1,1,\ldots,1]\in \R^{n+l+m}.$
  We compare our proposed GSS, RGSS-I and RGSS-II PGMRES methods with the GMRES method and PGMRES methods with BD, BD-MINRES, DS, RDF, and SS preconditioners, respectively. The BD preconditioner is used in its exact form as given in \cite[equation (6)]{Susanne2023}. The implementation of all the preconditioners is done as in Example \ref{exam1}. We consider the grids of sizes $16 \times 16$, $32 \times 32$, $64 \times 64$, and $128 \times 128.$ 
  
\vspace{2mm}
            \noindent\textbf{Parameter selection:} The parameter choices for the tested preconditioner are taken as in Example \ref{exam1}. In addition, for the proposed preconditioners, we 
            consider $\tau=0.0001.$ The optimal parameter \(\om\) (denoted as \(\om_{\text{exp}}\)) is experimentally determined by testing values in the interval \([10,30]\) with a step size of one, selecting the value that minimizes CPU time.

\begin{figure}[ht!]
				\centering
				\begin{subfigure}[b]{0.31\textwidth}
					\centering
					\includegraphics[width=\textwidth]{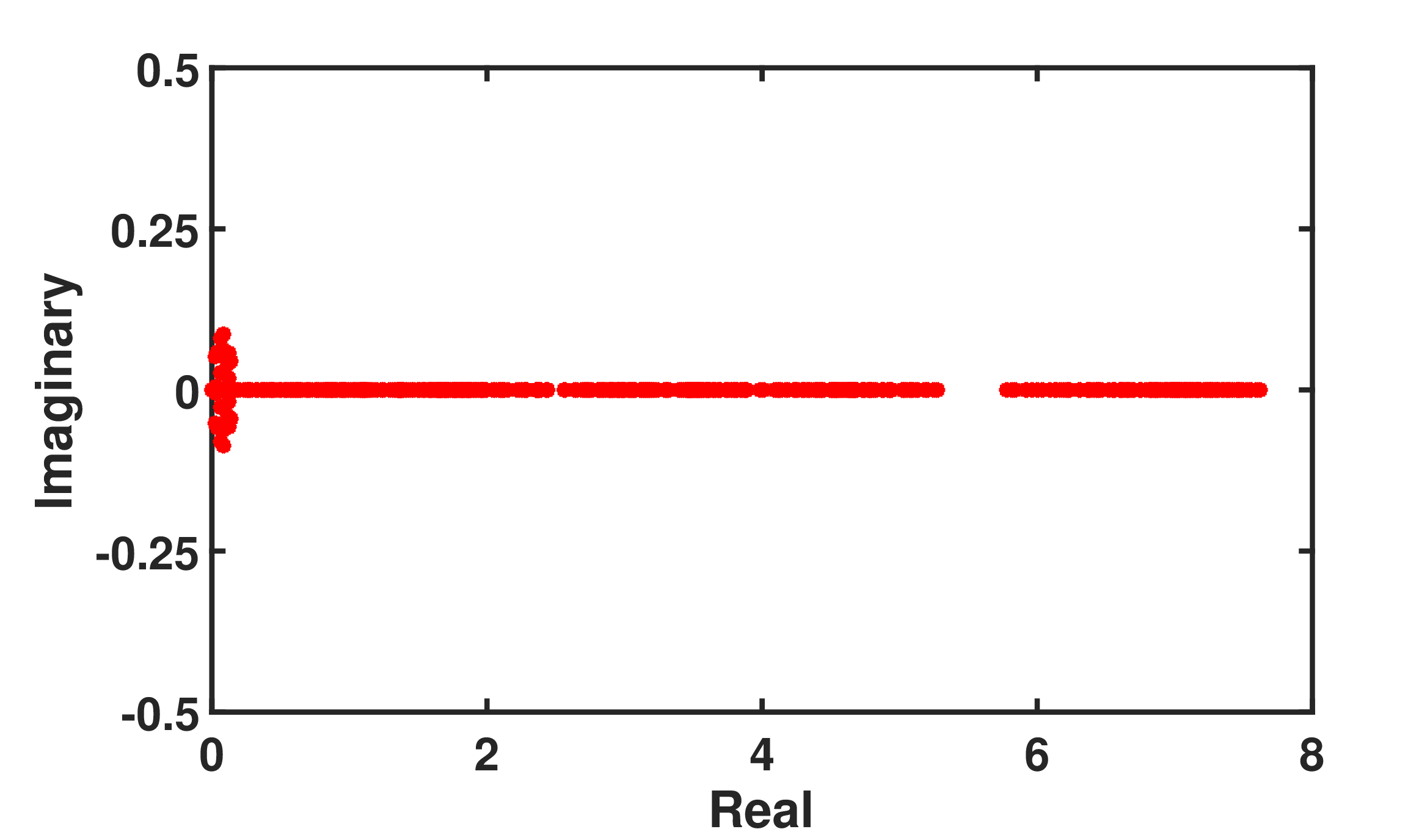}
					\caption{  $\B$ }
					\label{fig:original_L}
				\end{subfigure}
				\begin{subfigure}[b]{0.32\textwidth}
					\centering
					\includegraphics[width=\textwidth]{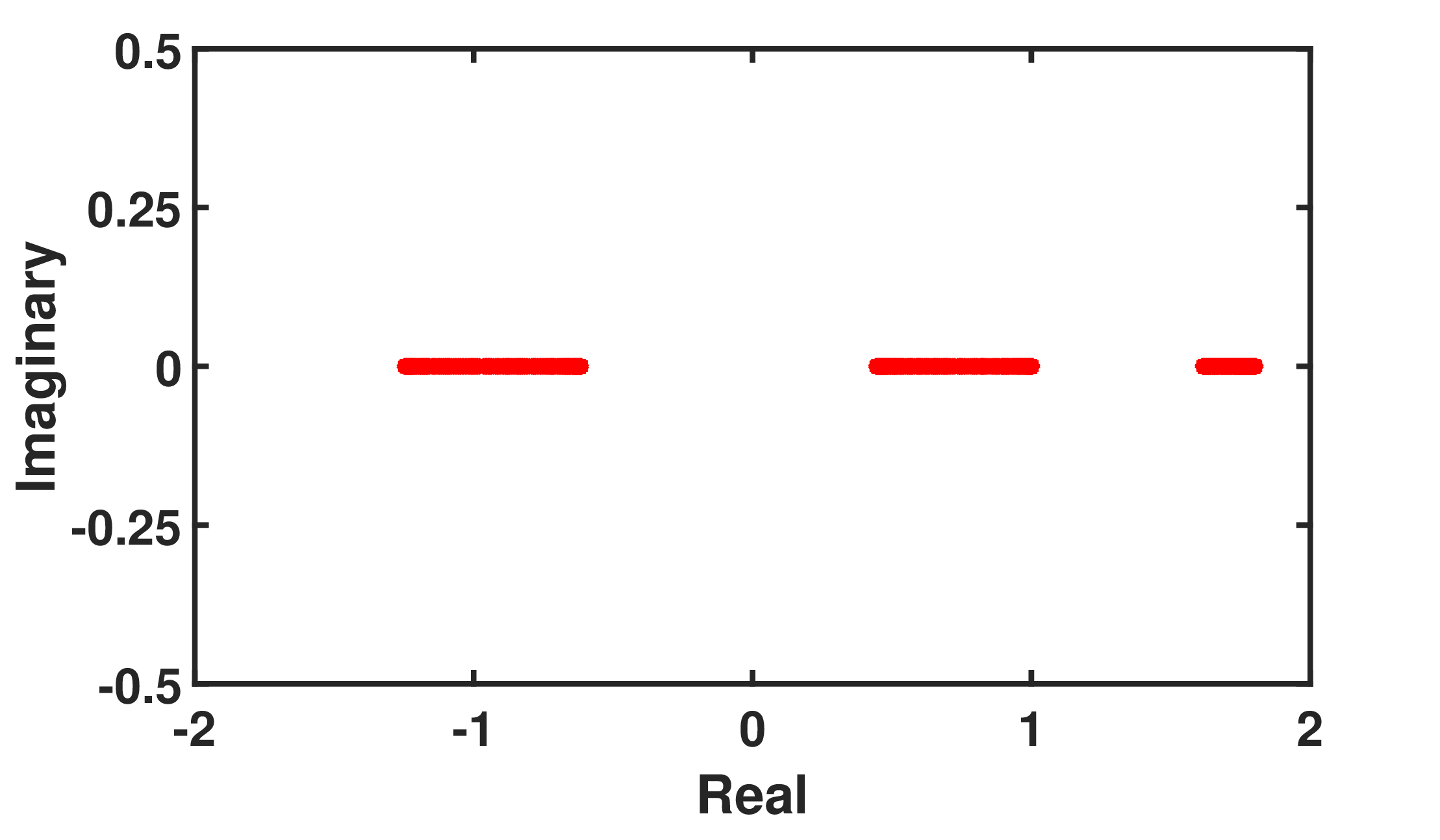}
					\caption{ $ \P_{\text{BD}}^{-1}\B$}
					\label{fig:BD_L}
				\end{subfigure}
                \begin{subfigure}[b]{0.31\textwidth}
					\centering
					\includegraphics[width=\textwidth]{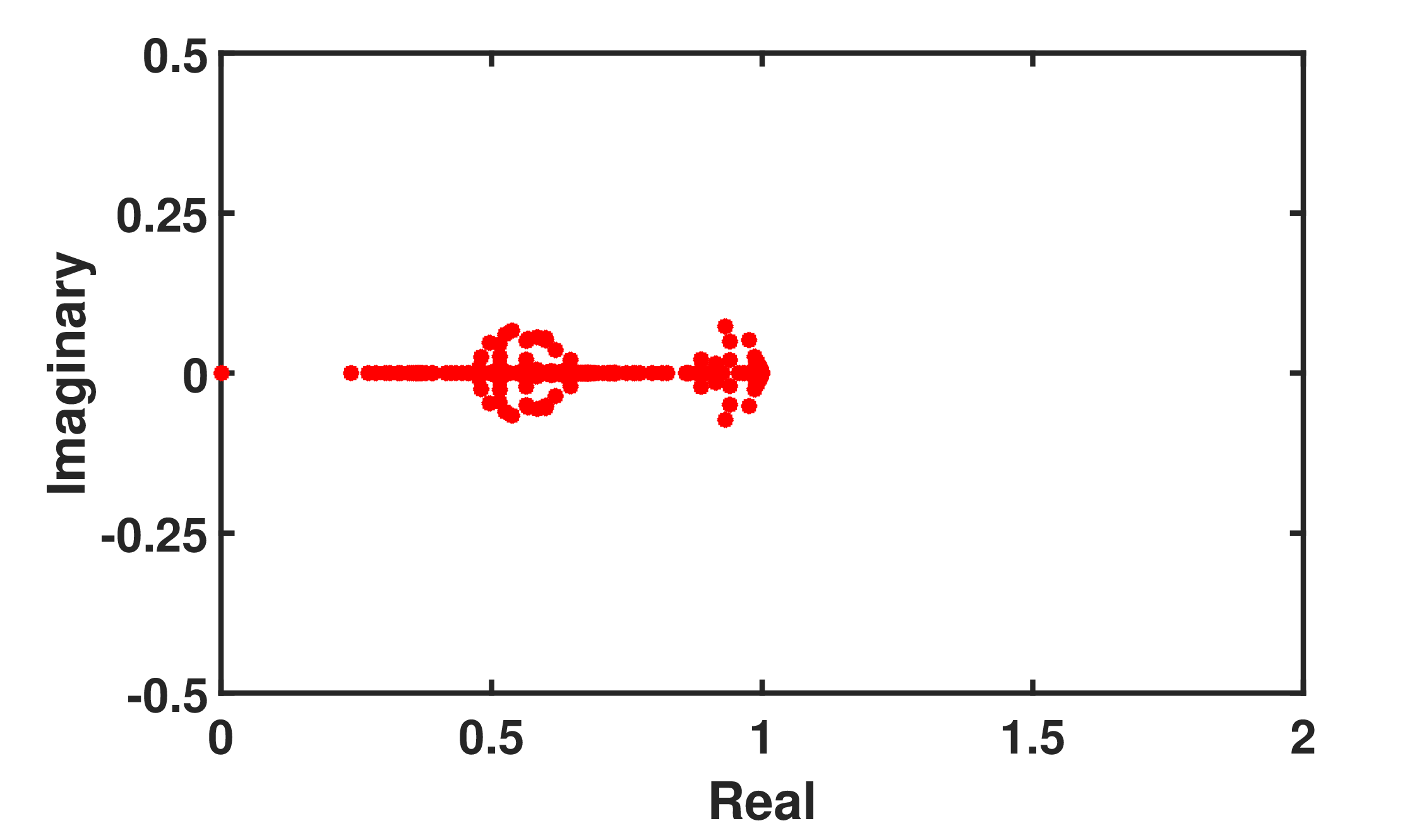}
					\caption{ $ \P_{\text{DS}}^{-1}\B$}
					\label{fig:DS_L}
				\end{subfigure}
    \begin{subfigure}[b]{0.28\textwidth}
					\centering
					\includegraphics[width=\textwidth]{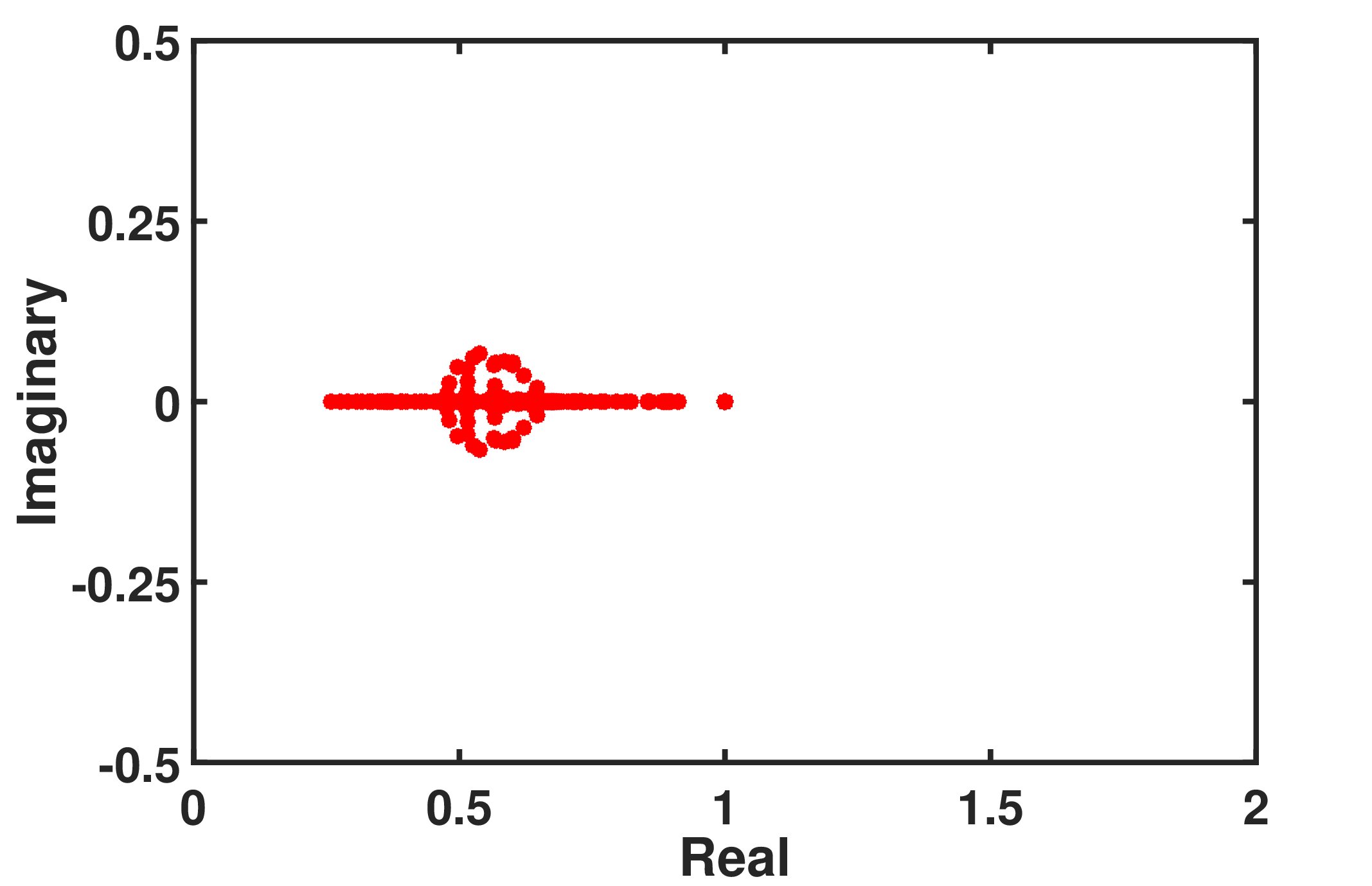}
					\caption{ $ \P_{\text{RDF}}^{-1}\B$}
					\label{fig:RDF_L}
				\end{subfigure}
    \begin{subfigure}[b]{0.31\textwidth}
					\centering
					\includegraphics[width=\textwidth]{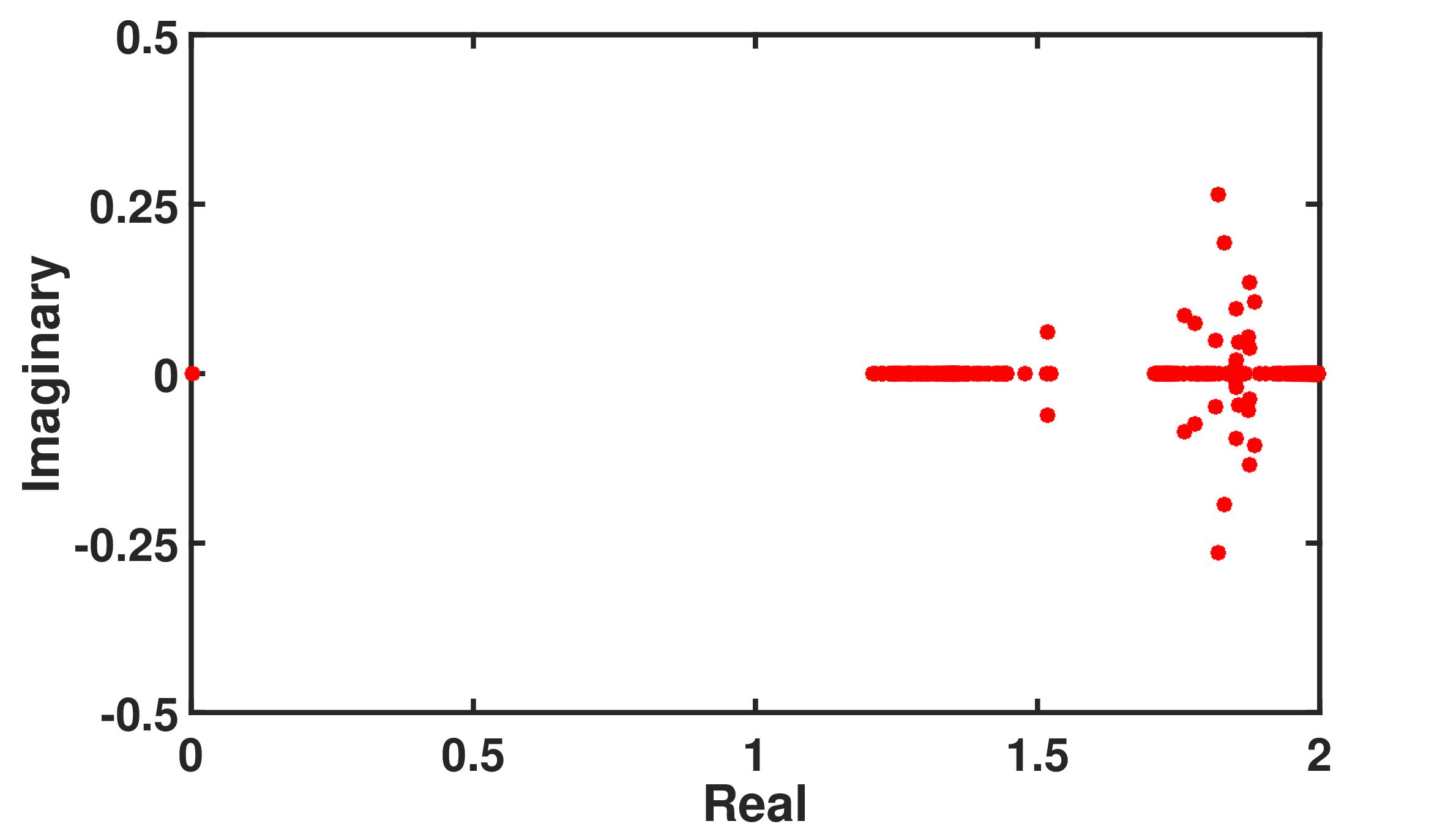}
					\caption{ $ \P_{\text{SS}}^{-1}\B$}
					\label{fig:SS_L}
				\end{subfigure}
    \begin{subfigure}[b]{0.33\textwidth}
					\centering
					\includegraphics[width=\textwidth]{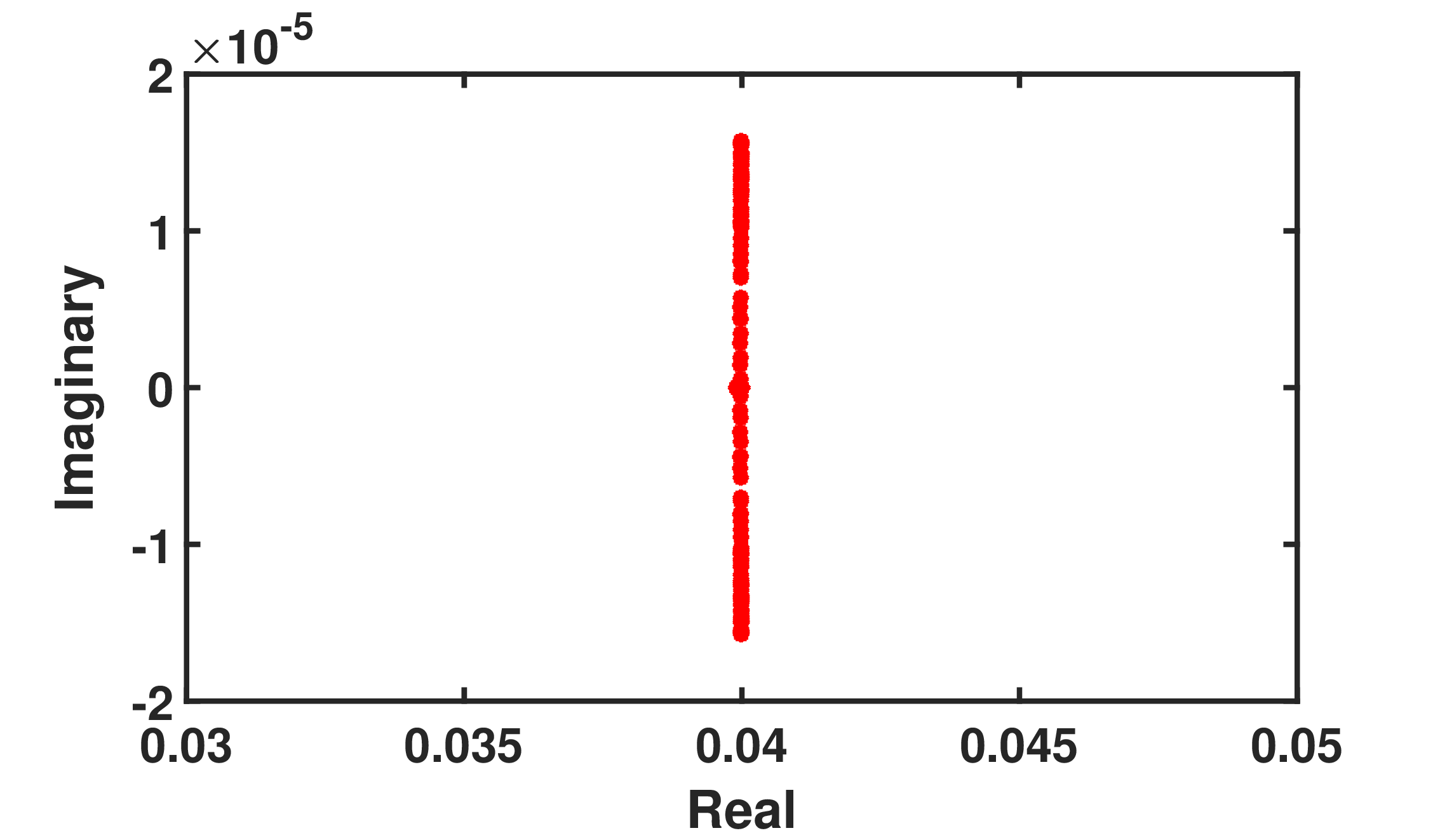}
					\caption{ $ \P_{\text{GSS}}^{-1}\B$   $(\om=25)$}
					\label{fig:GSS_L}
				\end{subfigure}
				\begin{subfigure}[b]{0.3\textwidth}
					\centering
					\includegraphics[width=\textwidth]{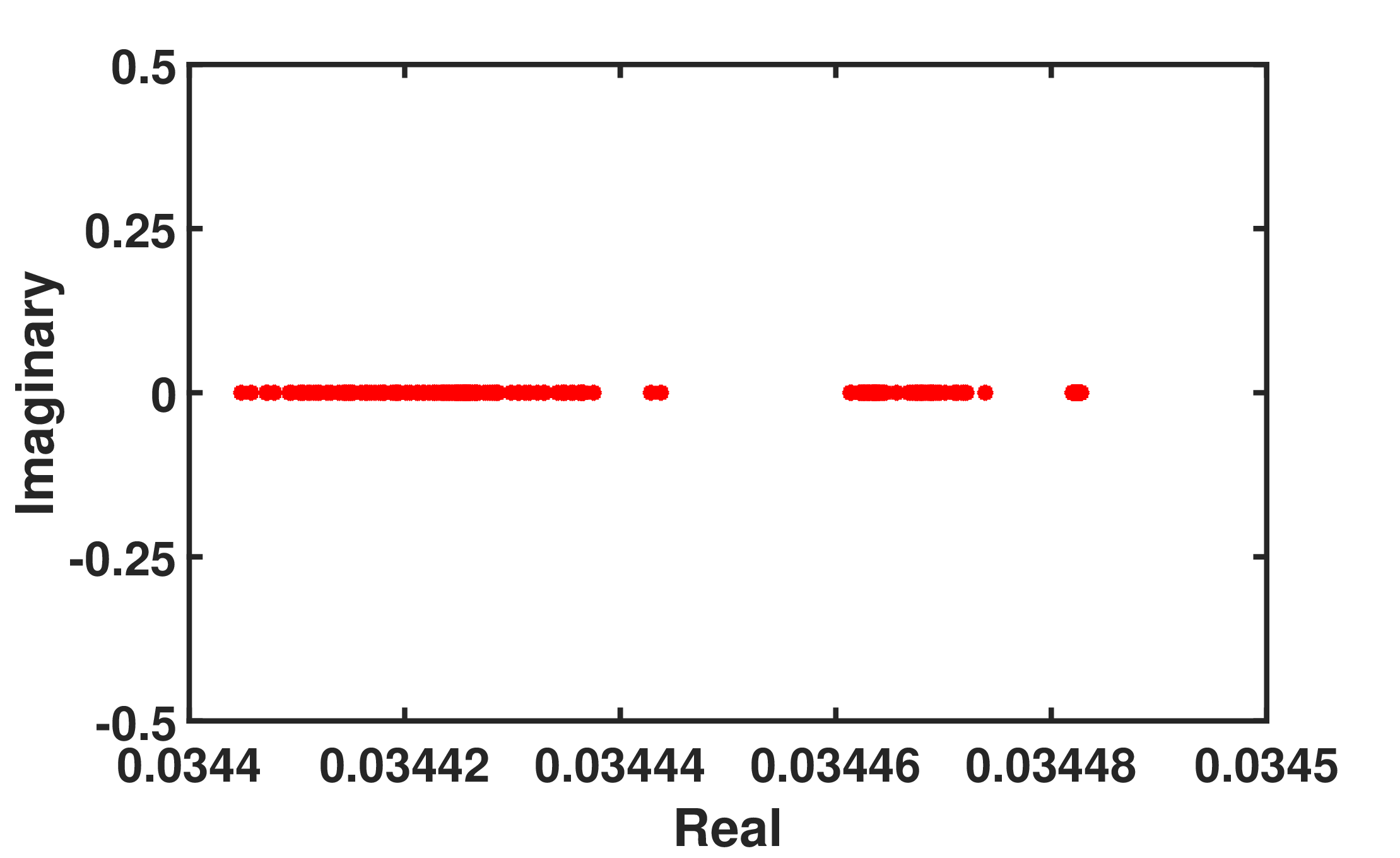}
					\caption{ $ \P_{\text{RGSS-I}}^{-1}\B$  $(\om=29)$ }
					\label{fig:RGSS_I_L}
				\end{subfigure}
    \begin{subfigure}[b]{0.33\textwidth}
					\centering
					\includegraphics[width=\textwidth]{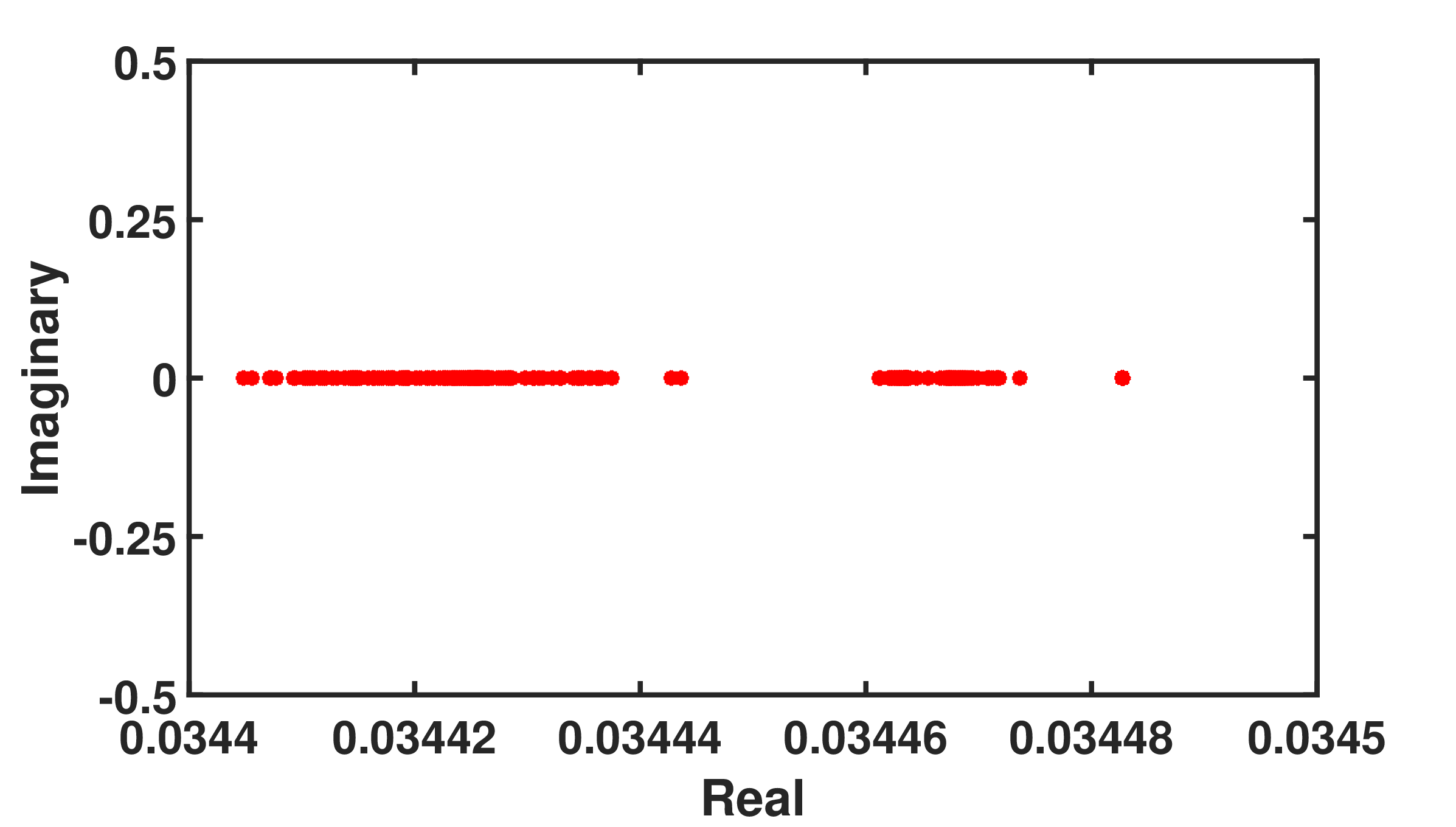}
					\caption{   $ \P_{\text{RGSS-II}}^{-1}\B$ $(\om=29)$ }
					\label{fig:RPSS_II_L}
				\end{subfigure}

				\caption{ Eigenvalue distributions of $\B,$ $ \P_{\mathrm{BD}}^{-1}\B,$ $ \P_{\mathrm{DS}}^{-1}\B,$ $  \P_{\mathrm{SS}}^{-1}\B,$ $ \P_{\mathrm{GSS}}^{-1}\B,$ $  \P_{\text{RGSS-I}}^{-1}\B$ and $ \P_{\text{RGSS-II}}^{-1}\B$  for Example \ref{EX2} with grid $16\times 16$.}
				\label{fig:Ex2}
			\end{figure}
            
            \vspace{2mm}
\noindent \textbf{Numerical results and discussions:} Table \ref{tab3} presents the numerical results for GMRES and various PGMRES methods. The standard GMRES method converges significantly slower than all PGMRES variants, requiring a substantially larger number of IT. In contrast, our proposed preconditioners outperform the compared ones in IT and CPU time. Notably, they achieve convergence within just 2 or 3 iterations, highlighting the efficiency of the GSS, RGSS-I, and RGSS-II preconditioners compared to the baseline preconditioners.
            
\vspace{2mm}
\noindent \textbf{Eigenvalue distributions:} To better illustrate the effectiveness of the proposed  preconditioners, Figure \ref{fig:Ex2} presents the eigenvalue distributions of $\B$ and preconditioned matrices $ \P_{\mathrm{BD}}^{-1}\B,$ $ \P_{\mathrm{DS}}^{-1}\B,$ $ \P_{\mathrm{RDF}}^{-1}\B,$ $ \P_{\mathrm{SS}}^{-1}\B,$ $ \P_{\text{GSS}}^{-1}\B\, {(\om=25)},$ $ \P_{\text{RGSS-I}}^{-1}\B\, {(\om=29)}$ and $ \P_{\text{RGSS-II}}^{-1}\B\, {(\om=29)}$ with ${\bm \nu}=0.1$ are displayed in Figure \ref{fig:Ex2}. Moreover, we computed the spectral radius of $\B$ and each preconditioned matrices, which are $\vartheta(\B)= 7.6198,$ $\vartheta(\P_{\mathrm{BD}}^{-1}\B)= 1.8018,$ $\vartheta(\P_{\mathrm{DS}}^{-1}\B)= 0.9997,$ $\vartheta(\P_{\mathrm{RDF}}^{-1}\B)= 1,$ $\vartheta(\P_{\mathrm{SS}}^{-1}\B)= 1.9974$, $\vartheta(\P_{\GSS}^{-1}\B)=0.04$, $\vartheta(\P_{\text{RGSS-I}}^{-1}\B)= 0.0345$, and $\vartheta(\P_{\text{RGSS-I}}^{-1}\B)=0.0345.$ Therefore, based on the spectral radius values and observation from Figure \ref{fig:Ex2}, it leads to the conclusion that the GSS, RGSS-I, and RGSS-II preconditioned matrices exhibit a more tightly clustered spectrum compared to other preconditioned matrices. This highlights the enhanced effectiveness of our proposed preconditioners.
 \end{exam}}
  
  \section{Concluding remarks}\label{sec:conclusion}
This paper introduces three preconditioners, termed GSS, RGSS-I and RGSS-II for solving DSPPs that arise in various applications. We provide a convergence analysis for the GSS iterative method, demonstrating that the method converges for any initial guess vector when the parameter $\om \geq 1/2.$ Moreover, spectral bounds for the preconditioned matrices are derived. Additionally, we have shown that the RGSS-II preconditioner requires at most $m+1$ iterations to solve the DSPP. Numerical experiments for the DSPP arising from the PDE-constrained optimization problem  {and leaky lid driven cavity problem} are performed, which demonstrate that the proposed preconditioners are efficient and outperform the existing state-of-the-art preconditioners.  {A line for future research direction is to investigate the spectral bounds of the RGSS-I and RGSS-II preconditioners when the diagonal blocks are nonsymmetric.}

  \section*{Acknowledgments}
	 {The authors sincerely acknowledge the anonymous reviewers for their constructive feedback to improve the paper}.	Pinki Khatun gratefully acknowledges the Council of Scientific $\&$ Industrial Research (CSIR), New Delhi, India, for providing financial support through a fellowship (File No. 09/1022(0098)/2020-EMR-I).
		
		\bibliography{reference}

\begin{thebibliography}{40}
\providecommand{\natexlab}[1]{#1}
\providecommand{\url}[1]{\texttt{#1}}
\expandafter\ifx\csname urlstyle\endcsname\relax
  \providecommand{\doi}[1]{doi: #1}\else
  \providecommand{\doi}{doi: \begingroup \urlstyle{rm}\Url}\fi

\bibitem[Abdolmaleki et~al.(2022)Abdolmaleki, Karimi, and Salkuyeh]{BDMaryam}
M.~Abdolmaleki, S.~Karimi, and D.~K. Salkuyeh.
\newblock A new block-diagonal preconditioner for a class of {$3\times 3$}
  block saddle point problems.
\newblock \emph{Mediterr. J. Math.}, 19\penalty0 (1):\penalty0 43, 15, 2022.

\bibitem[Ahmad and Khatun(2025)]{Pinki_PESS}
S.~S. Ahmad and P.~Khatun.
\newblock A robust parameterized enhanced shift-splitting preconditioner for
  three-by-three block saddle point problems.
\newblock \emph{Journal of Computational and Applied Mathematics},
  459:\penalty0 116358, 2025.

\bibitem[Bai(2019)]{HSS2019}
Z.-Z. Bai.
\newblock Regularized {HSS} iteration methods for stabilized saddle-point
  problems.
\newblock \emph{IMA J. Numer. Anal.}, 39\penalty0 (4):\penalty0 1888--1923,
  2019.

\bibitem[Bai and Wang(2008)]{ZZBaiUzawa}
Z.-Z. Bai and Z.-Q. Wang.
\newblock Some preconditiong inexact {U}zawa methods for generalized saddle
  point problems.
\newblock \emph{Linear Algebra Appl.}, 428\penalty0 (11-12):\penalty0
  2900--2932, 2008.

\bibitem[Bai et~al.(2003)Bai, Golub, and Ng]{HSS}
Z.-Z. Bai, G.~H. Golub, and M.~K. Ng.
\newblock Hermitian and skew-{H}ermitian splitting methods for non-{H}ermitian
  positive definite linear systems.
\newblock \emph{SIAM J. Matrix Anal. Appl.}, 24\penalty0 (3):\penalty0
  603--626, 2003.

\bibitem[Bai et~al.(2005)Bai, Parlett, and Wang]{GSOR2005}
Z.-Z. Bai, B.~N. Parlett, and Z.-Q. Wang.
\newblock On generalized successive overrelaxation methods for augmented linear
  systems.
\newblock \emph{Numer. Math.}, 102\penalty0 (1):\penalty0 1--38, 2005.

\bibitem[Bai et~al.(2006)Bai, Yin, and Su]{BaiSS}
Z.-Z. Bai, J.-F. Yin, and Y.-F. Su.
\newblock A shift-splitting preconditioner for non-{H}ermitian positive
  definite matrices.
\newblock \emph{J. Comput. Math.}, 24\penalty0 (4):\penalty0 539--552, 2006.

\bibitem[Benzi and Guo(2011)]{BenziDS2011}
M.~Benzi and X.-P. Guo.
\newblock A dimensional split preconditioner for {S}tokes and linearized
  {N}avier-{S}tokes equations.
\newblock \emph{Appl. Numer. Math.}, 61\penalty0 (1):\penalty0 66--76, 2011.

\bibitem[Benzi and Simoncini(2006)]{Benzi2006}
M.~Benzi and V.~Simoncini.
\newblock On the eigenvalues of a class of saddle point matrices.
\newblock \emph{Numer. Math.}, 103\penalty0 (2):\penalty0 173--196, 2006.

\bibitem[Benzi et~al.(2005)Benzi, Golub, and Liesen]{Benzi2005}
M.~Benzi, G.~H. Golub, and J.~Liesen.
\newblock Numerical solution of saddle point problems.
\newblock \emph{Acta Numer.}, 14:\penalty0 1--137, 2005.

\bibitem[Benzi et~al.(2011)Benzi, Ng, Niu, and Wang]{BenziRDF2011}
M.~Benzi, M.~Ng, Q.~Niu, and Z.~Wang.
\newblock A relaxed dimensional factorization preconditioner for the
  incompressible {N}avier-{S}tokes equations.
\newblock \emph{J. Comput. Phys.}, 230\penalty0 (16):\penalty0 6185--6202,
  2011.

\bibitem[Bojanczyk et~al.(2003)Bojanczyk, Higham, and Patel]{ILSE2003}
A.~Bojanczyk, N.~J. Higham, and H.~Patel.
\newblock The equality constrained indefinite least squares problem: theory and
  algorithms.
\newblock \emph{BIT}, 43\penalty0 (3):\penalty0 505--517, 2003.

\bibitem[Bradley and Greif(2023)]{Susanne2023}
S.~Bradley and C.~Greif.
\newblock Eigenvalue bounds for double saddle-point systems.
\newblock \emph{IMA J. Numer. Anal.}, 43\penalty0 (6):\penalty0 3564--3592,
  2023.

\bibitem[Cao(2019)]{CAOSS19}
Y.~Cao.
\newblock Shift-splitting preconditioners for a class of block three-by-three
  saddle point problems.
\newblock \emph{Appl. Math. Lett.}, 96:\penalty0 40--46, 2019.

\bibitem[Cao et~al.(2014)Cao, Du, and Niu]{CaoSS2014}
Y.~Cao, J.~Du, and Q.~Niu.
\newblock Shift-splitting preconditioners for saddle point problems.
\newblock \emph{J. Comput. Appl. Math.}, 272:\penalty0 239--250, 2014.

\bibitem[Cao et~al.(2017)Cao, Miao, and Ren]{CaoSS2017}
Y.~Cao, S.-X. Miao, and Z.-R. Ren.
\newblock On preconditioned generalized shift-splitting iteration methods for
  saddle point problems.
\newblock \emph{Comput. Math. Appl.}, 74\penalty0 (4):\penalty0 859--872, 2017.

\bibitem[Elman et~al.(2007)Elman, Ramage, and Silvester]{Matlab2007}
H.~C. Elman, A.~Ramage, and D.~J. Silvester.
\newblock Algorithm 886: {IFISS}, a {M}atlab toolbox for modelling
  incompressible flow.
\newblock \emph{ACM Trans. Math. Software}, 33\penalty0 (2):\penalty0 Art. 14,
  18, 2007.

\bibitem[Elman et~al.(2014)Elman, Silvester, and Wathen]{CFD2005}
H.~C. Elman, D.~J. Silvester, and A.~J. Wathen.
\newblock \emph{Finite Elements and Fast Iterative Solvers: with Applications
  in Incompressible Fluid Dynamics}.
\newblock Oxford University Press, Oxford, second edition, 2014.

\bibitem[Fan et~al.(2024)Fan, Li, Zhang, and Zhu]{SSPDE2024}
H.~Fan, Y.~Li, H.~Zhang, and X.~Zhu.
\newblock Preconditioners based on matrix splitting for the structured systems
  from elliptic {PDE}-constrained optimization problems.
\newblock \emph{Appl. Math. Comput.}, 463:\penalty0 Paper No. 128341, 8, 2024.

\bibitem[Golub et~al.(2001)Golub, Wu, and Yuan]{SOR2001}
G.~H. Golub, X.~Wu, and J.-Y. Yuan.
\newblock S{OR}-like methods for augmented systems.
\newblock \emph{BIT}, 41\penalty0 (1):\penalty0 71--85, 2001.

\bibitem[Grigori et~al.(2019)Grigori, Niu, and Xu]{SDS2019}
L.~Grigori, Q.~Niu, and Y.~Xu.
\newblock Stabilized dimensional factorization preconditioner for solving
  incompressible {N}avier-{S}tokes equations.
\newblock \emph{Appl. Numer. Math.}, 146:\penalty0 309--327, 2019.

\bibitem[Han and Yuan(2013)]{OPTMIZATION1}
D.~Han and X.~Yuan.
\newblock Local linear convergence of the alternating direction method of
  multipliers for quadratic programs.
\newblock \emph{SIAM J. Numer. Anal.}, 51\penalty0 (6):\penalty0 3446--3457,
  2013.

\bibitem[Huang(2020)]{NaHuangVUZ}
N.~Huang.
\newblock Variable parameter {U}zawa method for solving a class of block
  three-by-three saddle point problems.
\newblock \emph{Numer. Algorithms}, 85\penalty0 (4):\penalty0 1233--1254, 2020.

\bibitem[Huang and Ma(2019)]{HuangNA}
N.~Huang and C.-F. Ma.
\newblock Spectral analysis of the preconditioned system for the {$3\times3$}
  block saddle point problem.
\newblock \emph{Numer. Algor.}, 81\penalty0 (2):\penalty0 421--444, 2019.

\bibitem[Huang et~al.(2019)Huang, Dai, and Hu]{DUM2019}
N.~Huang, Y.-H. Dai, and Q.~Hu.
\newblock Uzawa methods for a class of block three-by-three saddle-point
  problems.
\newblock \emph{Numer. Linear Algebra Appl.}, 26\penalty0 (6):\penalty0 e2265,
  26, 2019.

\bibitem[Ke and Ma(2018)]{ChangFeng2018}
Y.-F. Ke and C.-F. Ma.
\newblock Some preconditioners for elliptic {PDE}-constrained optimization
  problems.
\newblock \emph{Comput. Math. Appl.}, 75\penalty0 (8):\penalty0 2795--2813,
  2018.

\bibitem[Li et~al.(2014)Li, Wu, Yang, and Yuan]{MUzawa2014}
X.~Li, Y.-J. Wu, A.-L. Yang, and J.-Y. Yuan.
\newblock Modified accelerated parameterized inexact {U}zawa method for
  singular and nonsingular saddle point problems.
\newblock \emph{Appl. Math. Comput.}, 244:\penalty0 552--560, 2014.

\bibitem[Mirchi and Salkuyeh(2020)]{Salkuyeh2020}
H.~Mirchi and D.~K. Salkuyeh.
\newblock A new preconditioner for elliptic {PDE}-constrained optimization
  problems.
\newblock \emph{Numer. Algorithms}, 83\penalty0 (2):\penalty0 653--668, 2020.

\bibitem[Pearson and Wathen(2012)]{Pearson2012}
J.~W. Pearson and A.~J. Wathen.
\newblock A new approximation of the {S}chur complement in preconditioners for
  {PDE}-constrained optimization.
\newblock \emph{Numer. Linear Algebra Appl.}, 19\penalty0 (5):\penalty0
  816--829, 2012.

\bibitem[Pearson et~al.(2012)Pearson, Stoll, and Wathen]{Pearson2}
J.~W. Pearson, M.~Stoll, and A.~J. Wathen.
\newblock Regularization-robust preconditioners for time-dependent
  {PDE}-constrained optimization problems.
\newblock \emph{SIAM J. Matrix Anal. Appl.}, 33\penalty0 (4):\penalty0
  1126--1152, 2012.

\bibitem[Pearson et~al.(2014)Pearson, Stoll, and Wathen]{Pearson3}
J.~W. Pearson, M.~Stoll, and A.~J. Wathen.
\newblock Preconditioners for state-constrained optimal control problems with
  {M}oreau-{Y}osida penalty function.
\newblock \emph{Numer. Linear Algebra Appl.}, 21\penalty0 (1):\penalty0 81--97,
  2014.

\bibitem[Rees()]{github}
T.~Rees.
\newblock Github - tyronerees/poisson-control. {A}ccessed: 2022-3-21, 2010.
\newblock \url{https://github.com/ tyronerees/poisson-control}.

\bibitem[Rees and Stoll(2010)]{ReesBT2010}
T.~Rees and M.~Stoll.
\newblock Block-triangular preconditioners for {PDE}-constrained optimization.
\newblock \emph{Numer. Linear Algebra Appl.}, 17\penalty0 (6):\penalty0
  977--996, 2010.

\bibitem[Rees et~al.(2010)Rees, Dollar, and Wathen]{PDE-constrained2010}
T.~Rees, H.~S. Dollar, and A.~J. Wathen.
\newblock Optimal solvers for {PDE}-constrained optimization.
\newblock \emph{SIAM J. Sci. Comput.}, 32\penalty0 (1):\penalty0 271--298,
  2010.

\bibitem[Saad(2003)]{YSAAD}
Y.~Saad.
\newblock \emph{Iterative Methods for Sparse Linear Systems}.
\newblock SIAM, Philadelphia, PA, second edition, 2003.

\bibitem[Saad and Schultz(1986)]{gmres}
Y.~Saad and M.~H. Schultz.
\newblock G{MRES}: a generalized minimal residual algorithm for solving
  nonsymmetric linear systems.
\newblock \emph{SIAM J. Sci. Statist. Comput.}, 7\penalty0 (3):\penalty0
  856--869, 1986.

\bibitem[Salkuyeh et~al.(2015)Salkuyeh, Masoudi, and Hezari]{GSS2015}
D.~K. Salkuyeh, M.~Masoudi, and D.~Hezari.
\newblock On the generalized shift-splitting preconditioner for saddle point
  problems.
\newblock \emph{Appl. Math. Lett.}, 48:\penalty0 55--61, 2015.

\bibitem[Yang et~al.(2024)Yang, Zhu, and Yu-Jiang]{MDS2024}
A.-L. Yang, J.-L. Zhu, and W.~Yu-Jiang.
\newblock Multi-parameter dimensional split preconditioner for three-by-three
  block system of linear equations.
\newblock \emph{Numer. Algorithms}, 95\penalty0 (2):\penalty0 721--745, 2024.

\bibitem[Zhang et~al.(2014)Zhang, Yang, and Wang]{Uzawa2014}
G.-F. Zhang, J.-L. Yang, and S.-S. Wang.
\newblock On generalized parameterized inexact {U}zawa method for a block
  two-by-two linear system.
\newblock \emph{J. Comput. Appl. Math.}, 255:\penalty0 193--207, 2014.

\bibitem[Zhang and Huang(2014)]{BCD2014}
X.~Zhang and Y.~Huang.
\newblock On block preconditioners for {PDE}-constrained optimization problems.
\newblock \emph{J. Comput. Math.}, 32\penalty0 (3):\penalty0 272--283, 2014.

\end{thebibliography}
		\bibliographystyle{abbrvnat}
		
	\end{document}